\documentclass[12pt,twoside]{article}

\usepackage{jltmac2e}

\renewcommand{\version}{{\rm Version of~}\date} 
\renewcommand{\version}{July 14, 2006}
\renewcommand{\vol}{??}                             
\renewcommand{\annum}{??}                           
\renewcommand{\title}{\centerline{Lie elements in pre-Lie algebras,
             trees and cohomology operations}}         
\renewcommand{\author}{M.~Markl\thanks{The author was supported by the grant GA \v CR 201/05/2117 and by
   the Academy of Sciences of the Czech Republic, 
   Institutional Research Plan No.~AV0Z10190503.}}                   
\renewcommand{\lastname}{Markl}    

\def\rec{} 
\def\fin{} 

\renewcommand{\address}{Author: M.~Markl

                       Mathematical Institute

                       Czech Academy of Sciences
 
                       \v Zitn\'a 25, 115 67 Praha 1

                       The Czech Republic

                        \catcode`\@=11
                        markl@math.cas.cz
                        \catcode`\@=13
}     

\begin{document} 
\firstpage


\newtheorem{odstavec}[Theorem]{\hskip -2mm}

\def\calP{{\cal P}}\def\calB{{\cal B}}\def\calS{{\cal S}}
\def\Aalg{\mbox {${\cal A}$-{\tt alg}}}
\def\Palg{\mbox {${\cal P}$-{\tt alg}}}
\def\pLie{p{\cal L{\it ie\/}}}\def\Brace{{\cal B{\it race\/}}}
\def\Prim{{\it Prim\/}}\def\calH{{\cal H}}
\def\Ass{{\cal A{\it ss\/}}} \def\Com{{\cal C{\it om\/}}}
\def\Lie{{\cal L{\it ie\/}}}\def\calC{{\cal C}}\def\calA{{\cal A}}
\def\qbezier{\bezier{100}}
\def\spin{
{
\unitlength=.017em
\begin{picture}(30.00,40.00)(0.00,0.00)
\qbezier(0.00,0.00)(20.00,0.00)(20.00,20.00)
\qbezier(0.00,40.00)(20.00,40.00)(20.00,20.00)
\qbezier(10.00,20.00)(10.00,0.00)(30.00,0.00)
\qbezier(30.00,40.00)(10.00,40.00)(10.00,20.00)
\end{picture}}\,
}
\def\ext{\mbox{\large$\land$}}
\def\Blie{{\cal B}_{\it Lie}}\def\Bshl{{\cal B}_{L_\infty}}
\def\Coder{{\it Coder\/}}\def\desusp{{\downarrow \hskip .1em}}
\def\CE#1{C_{\it CE}^{#1}(L;L)}\def\dce{{d_{\it CE}}}
\def\Rada#1#2#3{#1_{#2},\dots,#1_{#3}}
\def\Bl{{\it bl\/}}\def\Span{{\it Span\/}}
\def\brpL{{\overline {\rpL}}}
\def\Tr{{\rm Tr}}\def\otexp#1#2{#1^{\otimes #2}}
\def\oh{{\circ}}\def\td{{\tilde d}}\def\epi{\to \hskip -.6em \to}
\def\ccdot{\hskip .2em{\mbox {\small $\bullet$}}\hskip .2em}
\def\ot{\otimes}\def\Sh{{\it Sh}}\def\pl{\pL}
\def\T{{\sf T}}\def\fLie{{\sf L}}\def\fpLie{p{\sf L}}\def\bfk{{\bf k}}
\def\fL{\fLie}\def\fA{\T}\def\bDelta{{\overline{\Delta}}} 
\def\bfA{{\overline{\T}}}\def\Ker{{\it Ker}}\def\Im{{\it Im}}
\def\bar{{\hskip .2em\raisebox{-.2em}{\rule{.1em}{1em}}\hskip .2em}}
\def\L{\fLie}\def\bT{\bfA}\def\pL{{\sf pL}} \def\rpL{{\rm r\pL}}
\def\root{\rule{.4em}{.4em}}
\def\cases#1#2#3#4{
                  \left\{
                         \begin{array}{ll}
                           #1,\ &\mbox{#2}
                           \\
                           #3,\ &\mbox{#4}
                          \end{array}
                   \right.
}

\def\stromI#1#2{
\unitlength .5em
\thicklines
\ifx\plotpoint\undefined\newsavebox{\plotpoint}\fi 
\begin{picture}(7.03,5)(0,9)
\put(6.,13){\line(0,-1){3.97}}
\put(6.,13){\makebox(0,0)[cc]{$\root$}}
\put(6.,9){\makebox(0,0)[cc]{$\bullet$}}
\put(7.,13){\makebox(0,0)[lc]{$#1$}}
\put(7.,9){\makebox(0,0)[lc]{$#2$}}
\end{picture}
}

\def\stromII#1#2#3{
{
\thicklines
\unitlength=.95pt
\begin{picture}(10.00,40.00)(0.00,0.00)
\put(10.00,0.00){\makebox(0.00,0.00)[l]{$#3$}}
\put(10.00,20.00){\makebox(0.00,0.00)[l]{$#2$}}
\put(10.00,40.00){\makebox(0.00,0.00)[l]{$#1$}}
\put(0.00,0.00){\makebox(0.00,0.00){$\bullet$}}
\put(0.00,20.00){\makebox(0.00,0.00){$\bullet$}}
\put(0.00,40.00){\makebox(0.00,0.00){$\root$}}
\put(0.00,40.00){\line(0,-1){40.00}}
\end{picture}}
}

\def\stromIII#1#2#3{
{
\unitlength=1.000000pt
\thicklines
\begin{picture}(40.00,40.00)(0.00,0.00)
\put(40.00,10.00){\makebox(0.00,0.00){$\bullet$}}
\put(0.00,10.00){\makebox(0.00,0.00){$\bullet$}}
\put(20.00,30.00){\makebox(0.00,0.00){$\root$}}
\put(20.00,40.00){\makebox(0.00,0.00){$#1$}}
\put(40.00,0.00){\makebox(0.00,0.00){$#3$}}
\put(0.00,0.00){\makebox(0.00,0.00){$#2$}}
\put(20.00,30.00){\line(1,-1){20.00}}
\put(0.00,10.00){\line(1,1){20.00}}
\end{picture}}
}

\def\stromIV#1#2#3#4{
{
\unitlength=1.000000pt
\thicklines
\begin{picture}(20.00,30.00)(0.00,30.00)
\put(10.00,0.00){\makebox(0.00,0.00)[l]{$#4$}}
\put(10.00,20.00){\makebox(0.00,0.00)[l]{$#3$}}
\put(10.00,40.00){\makebox(0.00,0.00)[l]{$#2$}}
\put(10.00,60.00){\makebox(0.00,0.00)[l]{$#1$}}
\put(0.00,60.00){\makebox(0.00,0.00){$\root$}}
\put(0.00,40.00){\makebox(0.00,0.00){$\bullet$}}
\put(0.00,20.00){\makebox(0.00,0.00){$\bullet$}}
\put(0.00,0.00){\makebox(0.00,0.00){$\bullet$}}
\put(0.00,60.00){\line(0,-1){60.00}}
\end{picture}}
}

\def\stromV#1#2#3#4{
{
\unitlength=1.000000pt
\thicklines
\begin{picture}(50.00,30.00)(0.00,30.00)
\put(20.00,50.00){\makebox(0.00,0.00){$\root$}}
\put(40.00,10.00){\makebox(0.00,0.00){$\bullet$}}
\put(40.00,30.00){\makebox(0.00,0.00){$\bullet$}}
\put(0.00,30.00){\makebox(0.00,0.00){$\bullet$}}
\put(40.00,0.00){\makebox(0.00,0.00){$#4$}}
\put(50.00,30.00){\makebox(0.00,0.00)[l]{$#3$}}
\put(0.00,20.00){\makebox(0.00,0.00)[c]{$#2$}}
\put(20.00,60.00){\makebox(0.00,0.00)[c]{$#1$}}
\put(40.00,30.00){\line(0,-1){20.00}}
\put(20.00,50.00){\line(1,-1){20.00}}
\put(0.00,30.00){\line(1,1){20.00}}
\end{picture}}
}

\def\stromVI#1#2#3#4{
{
\unitlength=1.000000pt
\thicklines
\begin{picture}(40.00,30.00)(0.00,30.00)
\put(40.00,0.00){\makebox(0.00,0.00){$#4$}}
\put(0.00,0.00){\makebox(0.00,0.00){$#3$}}
\put(30.00,30.00){\makebox(0.00,0.00)[l]{$#2$}}
\put(20.00,60.00){\makebox(0.00,0.00){$#1$}}
\put(40.00,10.00){\makebox(0.00,0.00){$\bullet$}}
\put(0.00,10.00){\makebox(0.00,0.00){$\bullet$}}
\put(20.00,30.00){\makebox(0.00,0.00){$\bullet$}}
\put(20.00,50.00){\makebox(0.00,0.00){$\root$}}
\put(20.00,30.00){\line(1,-1){20.00}}
\put(20.00,30.00){\line(-1,-1){20.00}}
\put(20.00,50.00){\line(0,-1){20.00}}
\end{picture}}
}

\def\stromVII#1#2#3#4{
{
\unitlength=1.000000pt
\thicklines
\begin{picture}(40.00,20.00)(0.00,10.00)
\put(40.00,10.00){\makebox(0.00,0.00){$\bullet$}}
\put(20.00,10.00){\makebox(0.00,0.00){$\bullet$}}
\put(0.00,10.00){\makebox(0.00,0.00){$\bullet$}}
\put(20.00,30.00){\makebox(0.00,0.00){$\root$}}
\put(20.00,40.00){\makebox(0.00,0.00){$#1$}}
\put(0.00,0.00){\makebox(0.00,0.00){$#2$}}
\put(20.00,0.00){\makebox(0.00,0.00){$#3$}}
\put(40.00,0.00){\makebox(0.00,0.00){$#4$}}
\put(20.00,30.00){\line(0,-1){20.00}}
\put(20.00,30.00){\line(1,-1){20.00}}
\put(0.00,10.00){\line(1,1){20.00}}
\end{picture}}
}

\def\stromVIII#1#2#3#4{
{
\unitlength=1.000000pt
\begin{picture}(10.00,20.00)(0.00,0.00)
\thicklines
\put(10.00,0.00){\makebox(0.00,0.00)[l]{$#4$}}
\put(10.00,20.00){\makebox(0.00,0.00)[l]{$#3$}}
\put(0.00,0.00){\makebox(0.00,0.00){$#2$}}
\put(0.00,20.00){\makebox(0.00,0.00){$#1$}}
\put(0.00,18.10){\line(0,-1){16.10}}
\end{picture}}
}

\def\stromIX#1#2#3#4#5#6{
{
\unitlength=.95pt
\begin{picture}(10.00,44.00)(0.00,0.00)
\thicklines
\put(10.00,0.00){\makebox(0.00,0.00)[l]{$#6$}}
\put(10.00,20.00){\makebox(0.00,0.00)[l]{$#5$}}
\put(10.00,40.00){\makebox(0.00,0.00)[l]{$#4$}}
\put(0.00,0.00){\makebox(0.00,0.00){$#3$}}
\put(0.00,20.00){\makebox(0.00,0.00){$#2$}}
\put(0.00,40.00){\makebox(0.00,0.00){$#1$}}
\put(0.00,18.00){\line(0,-1){16.00}}
\put(0.00,38.00){\line(0,-1){16.00}}
\end{picture}}
}

\def\stromX#1#2#3#4#5#6{
{
\unitlength=1.000000pt
\begin{picture}(40.00,40.00)(0.00,0.00)
\thicklines
\put(40.00,0.00){\makebox(0.00,0.00){$#6$}}
\put(20.00,40.00){\makebox(0.00,0.00){$#5$}}
\put(0.00,0.00){\makebox(0.00,0.00){$#4$}}
\put(40.00,10.00){\makebox(0.00,0.00){$#3$}}
\put(20.00,30.00){\makebox(0.00,0.00){$#2$}}
\put(0.00,10.00){\makebox(0.00,0.00){$#1$}}
\put(38.50,11.50){\line(-1,1){17.00}}
\put(1.50,11.50){\line(1,1){17.00}}
\end{picture}}
}

\def\stromXI#1#2#3#4#5#6{
{
\unitlength=1.000000pt
\begin{picture}(50.00,20.00)(-5.00,20.00)
\thicklines
\put(40.00,0.00){\makebox(0.00,0.00){$#6$}}
\put(20.00,40.00){\makebox(0.00,0.00){$#5$}}
\put(0.00,0.00){\makebox(0.00,0.00){$#4$}}
\put(40.00,10.00){\makebox(0.00,0.00){$#3$}}
\put(20.00,30.00){\makebox(0.00,0.00){$#2$}}
\put(0.00,10.00){\makebox(0.00,0.00){$#1$}}
\put(38.50,11.50){\line(-1,1){17}}
\put(0,11.50){\line(1,1){18.1}}
\put(1.5,10){\line(1,1){18.1}}
\end{picture}}
}

\def\stromXII
{
\unitlength 3mm
\begin{picture}(6.5,2.5)(9.8,13.5)
\thicklines
\put(15,12){\line(0,1){2}}
\put(11.2,14){\line(1,1){1.8}}
\put(10.8,14){\line(1,1){1.95}}
\put(15,14){\line(-1,1){1.85}}
\put(13,16){\makebox(0,0)[cc]{$\oh$}}
\put(11,14){\makebox(0,0)[cc]{$\bullet$}}
\put(15,14){\makebox(0,0)[cc]{$\bullet$}}
\put(15,12){\makebox(0,0)[cc]{$\bullet$}}
\put(10,14){\makebox(0,0)[cc]{$u$}}
\put(16,14){\makebox(0,0)[cc]{$v$}}
\put(16,12){\makebox(0,0)[cc]{$w$}}
\end{picture}
}

\def\stromXii#1#2#3
{
\unitlength 3mm
\begin{picture}(6.5,2.5)(9.8,13.5)
\thicklines
\put(15,12){\line(0,1){2}}
\put(11,14){\line(1,1){1.85}}
\put(15,14){\line(-1,1){1.85}}
\put(13,16){\makebox(0,0)[cc]{$\oh$}}
\put(11,14){\makebox(0,0)[cc]{$\bullet$}}
\put(15,14){\makebox(0,0)[cc]{$\bullet$}}
\put(15,12){\makebox(0,0)[cc]{$\bullet$}}
\put(10,14){\makebox(0,0)[cc]{$#1$}}
\put(16,14){\makebox(0,0)[cc]{$#2$}}
\put(16,12){\makebox(0,0)[cc]{$#3$}}
\end{picture}
}

\def\stromXIII
{
\unitlength 3mm
\begin{picture}(5.5,2.5)(-3,0)
\thicklines
\put(0,2){\line(0,-1){1.7}}
\put(-1.8,-2){\line(1,1){1.8}}
\put(-2.2,-2){\line(1,1){1.95}}
\put(2,-2){\line(-1,1){1.85}}
\put(0,2){\makebox(0,0)[cc]{$\bullet$}}
\put(-2,-2){\makebox(0,0)[cc]{$\bullet$}}
\put(2,-2){\makebox(0,0)[cc]{$\bullet$}}
\put(0,0){\makebox(0,0)[cc]{$\oh$}}
\put(0,3){\makebox(0,0)[cc]{$u$}}
\put(-2,-3){\makebox(0,0)[cc]{$v$}}
\put(2,-3){\makebox(0,0)[cc]{$w$}}
\end{picture}
}

\def\stromXIV
{
\unitlength 3mm
\begin{picture}(6.5,2.5)(9.8,13.5)
\thicklines
\put(15,12){\line(0,1){2}}
\put(11,14){\line(1,1){1.85}}
\put(14.8,14){\line(-1,1){1.8}}
\put(15.2,14){\line(-1,1){1.95}}
\put(13,16){\makebox(0,0)[cc]{$\oh$}}
\put(11,14){\makebox(0,0)[cc]{$\bullet$}}
\put(15,14){\makebox(0,0)[cc]{$\bullet$}}
\put(15,12){\makebox(0,0)[cc]{$\bullet$}}
\put(10,14){\makebox(0,0)[cc]{$v$}}
\put(16,14){\makebox(0,0)[cc]{$u$}}
\put(16,12){\makebox(0,0)[cc]{$w$}}
\end{picture}
}

\def\stromXIVe
{
\unitlength 3mm
\begin{picture}(6.5,2.5)(9.8,13.5)
\thicklines
\put(15,12){\line(0,1){2}}
\put(11,14){\line(1,1){1.85}}
\put(14.8,14){\line(-1,1){1.8}}
\put(15.2,14){\line(-1,1){1.95}}
\put(13,16){\makebox(0,0)[cc]{$\oh$}}
\put(11,14){\makebox(0,0)[cc]{$\bullet$}}
\put(15,14){\makebox(0,0)[cc]{$\bullet$}}
\put(15,12){\makebox(0,0)[cc]{$\bullet$}}
\put(10,14){\makebox(0,0)[cc]{$u$}}
\put(16,14){\makebox(0,0)[cc]{$v$}}
\put(16,12){\makebox(0,0)[cc]{$w$}}
\put(14.5,15.5){\makebox(0,0)[lb]{$e$}}
\end{picture}
}

\def\stromXV
{
\unitlength 3mm
\begin{picture}(6.5,2.5)(9.8,13.5)
\thicklines
\put(15,12){\line(0,1){2}}
\put(11,14){\line(1,1){1.85}}
\put(14.8,14){\line(-1,1){1.8}}
\put(15.2,14){\line(-1,1){1.95}}
\put(13,16){\makebox(0,0)[cc]{$\oh$}}
\put(11,14){\makebox(0,0)[cc]{$\bullet$}}
\put(15,14){\makebox(0,0)[cc]{$\bullet$}}
\put(15,12){\makebox(0,0)[cc]{$\bullet$}}
\put(10,14){\makebox(0,0)[cc]{$w$}}
\put(16,14){\makebox(0,0)[cc]{$u$}}
\put(16,12){\makebox(0,0)[cc]{$v$}}
\end{picture}
}

\def\stromXVI
{
\unitlength 3mm
\begin{picture}(5.5,2.5)(-3,0)
\thicklines
\put(0.15,2){\line(0,-1){1.8}}
\put(-.15,2){\line(0,-1){1.8}}
\put(-2,-2){\line(1,1){1.85}}
\put(2,-2){\line(-1,1){1.85}}
\put(0,2){\makebox(0,0)[cc]{$\bullet$}}
\put(-2,-2){\makebox(0,0)[cc]{$\bullet$}}
\put(2,-2){\makebox(0,0)[cc]{$\bullet$}}
\put(0,0){\makebox(0,0)[cc]{$\oh$}}
\put(0,3){\makebox(0,0)[cc]{$u$}}
\put(-2,-3){\makebox(0,0)[cc]{$v$}}
\put(2,-3){\makebox(0,0)[cc]{$w$}}
\end{picture}
}

\def\stromeXtra
{
\unitlength 3mm
\begin{picture}(5.5,2.5)(-3,-2)
\thicklines
\put(0.15,-2){\line(0,1){1.8}}
\put(-.15,-2){\line(0,1){1.8}}
\put(-2,-2){\line(1,1){1.85}}
\put(2,-2){\line(-1,1){1.85}}
\put(0,-2){\makebox(0,0)[cc]{$\bullet$}}
\put(-2,-2){\makebox(0,0)[cc]{$\bullet$}}
\put(2,-2){\makebox(0,0)[cc]{$\bullet$}}
\put(0,0){\makebox(0,0)[cc]{$\oh$}}
\put(0,-3){\makebox(0,0)[cc]{$u$}}
\put(-2,-3){\makebox(0,0)[cc]{$v$}}
\put(2,-3){\makebox(0,0)[cc]{$w$}}
\end{picture}
}

\def\stromExtra
{
\unitlength 3mm
\begin{picture}(5.5,2.5)(-3,-1)
\thicklines
\put(0,-2){\line(0,1){1.8}}
\put(0,2){\makebox(0,0)[cc]{$\bullet$}}
\put(0,2.6){\makebox(0,0)[bc]{$t$}}
\put(0,0.3){\line(0,1){1.8}}
\put(-2,-2){\line(1,1){1.85}}
\put(2,-2){\line(-1,1){1.85}}
\put(0,-2){\makebox(0,0)[cc]{$\bullet$}}
\put(-2,-2){\makebox(0,0)[cc]{$\bullet$}}
\put(2,-2){\makebox(0,0)[cc]{$\bullet$}}
\put(0,0){\makebox(0,0)[cc]{$\oh$}}
\put(0,-3){\makebox(0,0)[cc]{$u$}}
\put(-2,-3){\makebox(0,0)[cc]{$v$}}
\put(2,-3){\makebox(0,0)[cc]{$w$}}
\end{picture}
}

\def\stromXvi#1#2#3
{
\unitlength 3mm
\begin{picture}(5.5,2.5)(-3,0)
\thicklines
\put(0,2){\line(0,-1){1.8}}
\put(-2,-2){\line(1,1){1.85}}
\put(2,-2){\line(-1,1){1.85}}
\put(0,2){\makebox(0,0)[cc]{$\bullet$}}
\put(-2,-2){\makebox(0,0)[cc]{$\bullet$}}
\put(2,-2){\makebox(0,0)[cc]{$\bullet$}}
\put(0,0){\makebox(0,0)[cc]{$\oh$}}
\put(0,3){\makebox(0,0)[cc]{$#1$}}
\put(-2,-3){\makebox(0,0)[cc]{$#2$}}
\put(2,-3){\makebox(0,0)[cc]{$#3$}}
\end{picture}
}

\def\stromXX
{
\unitlength 5mm
\begin{picture}(4,3)(-1,1)
\thicklines
\put(0,1){\line(1,1){.89}}
\put(2,1){\line(-1,1){.9}}
\put(2,0){\line(0,1){1}}
\put(1,2.15){\line(0,1){.9}}
\put(2,1){\makebox(0,0)[cc]{$\bullet$}}
\put(2,0){\makebox(0,0)[cc]{$\bullet$}}
\put(1,3){\makebox(0,0)[cc]{$\bullet$}}
\put(1.5,3){\makebox(0,0)[cc]{$t$}}
\put(0,1){\makebox(0,0)[cc]{$\bullet$}}
\put(-.5,1){\makebox(0,0)[cc]{$u$}}
\put(2.5,1){\makebox(0,0)[cc]{$v$}}
\put(2.5,0){\makebox(0,0)[cc]{$w$}}
\put(1,2){\makebox(0,0)[cc]{$\oh$}}
\end{picture}
}

\def\stromXXI
{
\unitlength 5mm
\begin{picture}(4,3)(-1,1)
\thicklines
\put(0,-1)
{
\put(0,1){\line(1,1){.88}}
\put(2,1){\line(-1,1){.9}}
\put(1,3){\line(0,1){1}}
\put(1,2.15){\line(0,1){.9}}
\put(2,1){\makebox(0,0)[cc]{$\bullet$}}
\put(1,3){\makebox(0,0)[cc]{$\bullet$}}
\put(1,4){\makebox(0,0)[cc]{$\bullet$}}
\put(1.5,3){\makebox(0,0)[cc]{$u$}}
\put(1.5,4){\makebox(0,0)[cc]{$t$}}
\put(0,1){\makebox(0,0)[cc]{$\bullet$}}
\put(-.5,1){\makebox(0,0)[cc]{$v$}}
\put(2.5,1){\makebox(0,0)[cc]{$w$}}
\put(1,2){\makebox(0,0)[cc]{$\oh$}}
}
\end{picture}
}

\def\stromXXII
{
\unitlength 5mm
\thicklines
\begin{picture}(4,3)(-1,1)
\put(0,1){
\put(0,1){\line(1,1){.88}}
\put(2,1){\line(-1,1){.9}}
\put(2,0){\line(0,1){1}}
\put(2,0){\line(0,-1){1}}
\put(2,1){\makebox(0,0)[cc]{$\bullet$}}
\put(2,0){\makebox(0,0)[cc]{$\bullet$}}
\put(2,-1){\makebox(0,0)[cc]{$\bullet$}}
\put(0,1){\makebox(0,0)[cc]{$\bullet$}}
\put(-.5,1){\makebox(0,0)[cc]{$t$}}
\put(2.5,1){\makebox(0,0)[cc]{$u$}}
\put(2.5,0){\makebox(0,0)[cc]{$v$}}
\put(2.5,-1){\makebox(0,0)[cc]{$w$}}
\put(1,2){\makebox(0,0)[cc]{$\oh$}}
}
\end{picture}
}

\def\stromXXIII
{
\unitlength 5mm
\begin{picture}(4,3)(-1,1)
\thicklines
\put(0,1){
\put(2,1){\line(-1,1){.88}}
\put(0,1){\line(1,1){.9}}
\put(2,0){\line(0,1){1}}
\put(2,0){\line(0,-1){1}}
\put(2,1){\makebox(0,0)[cc]{$\bullet$}}
\put(2,0){\makebox(0,0)[cc]{$\bullet$}}
\put(2,-1){\makebox(0,0)[cc]{$\bullet$}}
\put(0,1){\makebox(0,0)[cc]{$\bullet$}}
\put(-.5,1){\makebox(0,0)[cc]{$u$}}
\put(2.5,1){\makebox(0,0)[cc]{$t$}}
\put(2.5,0){\makebox(0,0)[cc]{$v$}}
\put(2.5,-1){\makebox(0,0)[cc]{$w$}}
\put(1,2){\makebox(0,0)[cc]{$\oh$}}
}
\end{picture}
}

\def\stromXXIV
{
\unitlength 5mm
\begin{picture}(4,3)(-1,1.5)
\thicklines
\put(0,-.5){
\put(0,1){\line(1,1){.88}}
\put(0,3){\line(1,1){1}}
\put(2,1){\line(-1,1){.9}}
\put(1,2.15){\line(0,1){1.9}}
\put(2,1){\makebox(0,0)[cc]{$\bullet$}}
\put(1,4){\makebox(0,0)[cc]{$\bullet$}}
\put(1.5,4){\makebox(0,0)[cc]{$t$}}
\put(0,1){\makebox(0,0)[cc]{$\bullet$}}
\put(0,3){\makebox(0,0)[cc]{$\bullet$}}
\put(-.5,1){\makebox(0,0)[cc]{$v$}}
\put(-.5,3){\makebox(0,0)[cc]{$u$}}
\put(2.5,1){\makebox(0,0)[cc]{$w$}}
\put(1,2){\makebox(0,0)[cc]{$\oh$}}
}
\end{picture}
}

\def\stromXXV
{
\unitlength 5mm
\begin{picture}(4,3)(-1,.5)
\thicklines
\put(0,.5){
\put(0,1){\line(1,1){.88}}
\put(2,1){\line(-1,1){.9}}
\put(2,-1){\line(0,1){2}}
\put(0,-1){\line(0,1){2}}
\put(2,1){\makebox(0,0)[cc]{$\bullet$}}
\put(0,-1){\makebox(0,0)[cc]{$\bullet$}}
\put(2,-1){\makebox(0,0)[cc]{$\bullet$}}
\put(0,1){\makebox(0,0)[cc]{$\bullet$}}
\put(-.5,1){\makebox(0,0)[cc]{$t$}}
\put(2.5,1){\makebox(0,0)[cc]{$v$}}
\put(2.5,-1){\makebox(0,0)[cc]{$w$}}
\put(-.5,-1){\makebox(0,0)[cc]{$u$}}
\put(1,2){\makebox(0,0)[cc]{$\oh$}}
}
\end{picture}
}

\def\stromXXVI
{
\unitlength 5mm
\begin{picture}(4,3)(-1,1)
\thicklines
\put(0,1){\line(1,1){.9}}
\put(2,1){\line(-1,1){.9}}
\put(2,0){\line(0,1){1}}
\put(1,2.13){\line(0,1){.9}}
\put(2,1){\makebox(0,0)[cc]{$\bullet$}}
\put(2,0){\makebox(0,0)[cc]{$\bullet$}}
\put(1,3){\makebox(0,0)[cc]{$\bullet$}}
\put(1.5,3){\makebox(0,0)[cc]{$t$}}
\put(0,1){\makebox(0,0)[cc]{$\bullet$}}
\put(-.5,1){\makebox(0,0)[cc]{$u$}}
\put(2.5,1){\makebox(0,0)[cc]{$v$}}
\put(2.5,0){\makebox(0,0)[cc]{$w$}}
\put(1,2){\makebox(0,0)[cc]{$\oh$}}
\end{picture}
}

\def\stromXXVIextra
{
\unitlength 5mm
\begin{picture}(4,3)(-1,.5)
\thicklines
\put(0,1){\line(1,1){.9}}
\put(2,1){\line(-1,1){.9}}
\put(2,0){\line(0,1){1}}
\put(1,1.9){\line(0,-1){.9}}
\put(2,1){\makebox(0,0)[cc]{$\bullet$}}
\put(2,0){\makebox(0,0)[cc]{$\bullet$}}
\put(1,1){\makebox(0,0)[cc]{$\bullet$}}
\put(1,.7){\makebox(0,0)[tc]{$t$}}
\put(0,1){\makebox(0,0)[cc]{$\bullet$}}
\put(-.5,1){\makebox(0,0)[cc]{$u$}}
\put(2.5,1){\makebox(0,0)[cc]{$v$}}
\put(2.5,0){\makebox(0,0)[cc]{$w$}}
\put(1,2){\makebox(0,0)[cc]{$\oh$}}
\end{picture}
}

\def\stromXXVII
{
\unitlength 5mm
\begin{picture}(4,3)(-1,1)
\thicklines
\put(0,1){\line(1,1){.9}}
\put(2,1){\line(-1,1){.9}}
\put(2,0){\line(0,1){1}}
\put(1,2.15){\line(0,1){.9}}
\put(2,1){\makebox(0,0)[cc]{$\bullet$}}
\put(2,0){\makebox(0,0)[cc]{$\bullet$}}
\put(1,3){\makebox(0,0)[cc]{$\bullet$}}
\put(1.5,3){\makebox(0,0)[cc]{$t$}}
\put(0,1){\makebox(0,0)[cc]{$\bullet$}}
\put(-.5,1){\makebox(0,0)[cc]{$v$}}
\put(2.5,1){\makebox(0,0)[cc]{$u$}}
\put(2.5,0){\makebox(0,0)[cc]{$w$}}
\put(1,2){\makebox(0,0)[cc]{$\oh$}}
\end{picture}
}

\def\stromXXVIII
{
\unitlength 5mm
\begin{picture}(4,3)(-1,1)
\thicklines
\put(0,1){\line(1,1){.9}}
\put(2,1){\line(-1,1){.87}}
\put(2,0){\line(0,1){1}}
\put(1,2.15){\line(0,1){.9}}
\put(2,1){\makebox(0,0)[cc]{$\bullet$}}
\put(2,0){\makebox(0,0)[cc]{$\bullet$}}
\put(1,3){\makebox(0,0)[cc]{$\bullet$}}
\put(1.5,3){\makebox(0,0)[cc]{$t$}}
\put(0,1){\makebox(0,0)[cc]{$\bullet$}}
\put(-.6,1){\makebox(0,0)[cc]{$w$}}
\put(2.5,1){\makebox(0,0)[cc]{$u$}}
\put(2.5,0){\makebox(0,0)[cc]{$v$}}
\put(1,2){\makebox(0,0)[cc]{$\oh$}}
\end{picture}
}

\def\stromXXIX
{
\unitlength 5mm
\begin{picture}(4,3)(-1,1)
\thicklines
\put(0,-1)
{
\put(0,1){\line(1,1){.9}}
\put(2,1){\line(-1,1){.9}}
\put(1,3){\line(0,1){1}}
\put(1,2.14){\line(0,1){.9}}
\put(2,1){\makebox(0,0)[cc]{$\bullet$}}
\put(1,3){\makebox(0,0)[cc]{$\bullet$}}
\put(1,4){\makebox(0,0)[cc]{$\bullet$}}
\put(1.5,3){\makebox(0,0)[cc]{$u$}}
\put(1.5,4){\makebox(0,0)[cc]{$t$}}
\put(0,1){\makebox(0,0)[cc]{$\bullet$}}
\put(-.5,1){\makebox(0,0)[cc]{$v$}}
\put(2.5,1){\makebox(0,0)[cc]{$w$}}
\put(1,2){\makebox(0,0)[cc]{$\oh$}}
}
\end{picture}
}

\def\stromXXX
{
\unitlength 5mm
\begin{picture}(4,3)(-1,1.5)
\thicklines
\put(0,-.5){
\put(0,3){\line(1,1){.87}}
\put(0,1){\line(1,1){1}}
\put(2,1){\line(-1,1){.9}}
\put(1,2){\line(0,1){1.85}}
\put(2,1){\makebox(0,0)[cc]{$\bullet$}}
\put(1,4){\makebox(0,0)[cc]{$\oh$}}
\put(0,1){\makebox(0,0)[cc]{$\bullet$}}
\put(0,3){\makebox(0,0)[cc]{$\bullet$}}
\put(-.5,1){\makebox(0,0)[cc]{$v$}}
\put(1.5,2){\makebox(0,0)[cc]{$u$}}
\put(-.5,3){\makebox(0,0)[cc]{$t$}}
\put(2.5,1){\makebox(0,0)[cc]{$w$}}
\put(1,2){\makebox(0,0)[cc]{$\bullet$}}
}
\end{picture}
}

\begin{abstract}
We give a simple characterization of Lie elements in free pre-Lie
algebras as elements of the kernel of a map between spaces of trees. 
We explain how this result is related to natural operations on the
Chevalley-Eilenberg complex of a Lie algebra. 
We also indicate a possible relation to Loday's theory of
triplettes.
\end{abstract}

\vbadness=100000
\section{Main results, motivations and generalizations}
\label{sec0}

All algebraic objects in this note will be defined over a field $\bfk$
of characteristic zero and $V$ will always denote a 
$\bfk$-vector space. We will sometimes use the formalism of operads
explained, for example,
in~\cite{markl-shnider-stasheff:book}. 
Sections~2,~3
and~4 containing the
main results, however, do not rely on this language.

Let $\pL(V)$ denote the free pre-Lie algebra generated by $V$ and
$\pL(V)_L$ the associated Lie algebra. We will focus on the Lie
algebra $\L(V) \subset \pL(V)$ generated in $\pL(V)_L$ by $V$, called
the {\em subalgebra of Lie elements\/} in $\pl(V)$.  It is
known~\cite{dzhu-lof:HHA02} that $\L(V)$ is (isomorphic to) the free
Lie algebra generated by $V$; we will give a new short proof of this
statement in Section~3. Our main result, Theorem~\ref{.},
describes $\L(V)$ as the kernel of a map 
\begin{equation}
\label{psano_v_Srni}
d : \pL(V) \to \pL^1(V),
\end{equation}
where $\pL^1(V)$ is the subspace of degree $+1$ elements in the free
graded pre-Lie algebra $\pL^*(V,\oh)$ generated by $V$ and a degree
$+1$ `dummy' variable $\oh$.
The map~(\ref{psano_v_Srni}) is later in the paper identified 
with a very simple map between
spaces of trees, see Proposition~\ref{abych_neonemocnel} and
Corollary~\ref{,}.

Theorem~\ref{.} and Corollary~\ref{,} have 
immediate applications to the analysis of natural operations on the
Chevalley-Eilenberg complex of a Lie algebra. In a future work we also
plan to prove that Theorem~\ref{.} implies that the only natural
multilinear operations on vector fields on smooth manifolds are, in stable
dimensions, iterations of the standard Jacobi bracket. There is
also a possible relation of the results of this paper with Loday's
theory of triplettes. In the rest of this introductory
section, we discuss some of these applications and 
motivations in more detail.

\begin{odstavec}
\label{mot}
{\rm
{\it Motivations.\/}
In~\cite{markl:de} we studied, among other things, the
differential graded (dg-) operad $\Blie^*$ of natural operations on
the Chevalley-Eilenberg complex of a Lie algebra with coefficients in
itself, along with its homotopy version 
$\Bshl^*$, the operad of natural operations on the
Chevalley-Eilenberg complex of an $L_\infty$-algebra (=~strongly
homotopy Lie algebra, see~\cite{lada-markl:CommAlg95}).
We proposed:

\begin{problem}
\label{-/}
Describe the homotopy types, in the non-abelian derived category,
of the dg-operads $\Blie^*$ and  $\Bshl^*$ of natural operations on the
Chevalley-Eilenberg complex.
\end{problem}

The following conjecture was proposed by D.~Tamarkin.

\begin{conjecture}
\label{Tamarkin}
The operad $\Blie^*$ has the homotopy type of the operad
$\Lie$ for Lie algebras.
\end{conjecture}

It turns out that $\Bshl^*$, which is tied to
$\Blie^*$ by the `forgetful' map $c : \Bshl^* \to \Blie^*$, contains a
dg-sub-operad $\rpL^* = (\rpL^*,d)$ generated by symmetric
braces~\cite{lm:sb} such that $\rpL^0$ (the sub-operad of degree $0$
elements) is the operad $\pLie$ governing pre-Lie algebras. Moreover,
both $\rpL^0$, $\rpL^1$ and the differential $d : \rpL^0 \to \rpL^1$
have very explicit descriptions in terms of planar trees.
Our conviction in Conjecture~\ref{Tamarkin} made us believe that the sub-operad
$$
H^0(\rpL^*) = \Ker \left(\rule{0em}{1em} 
d : \rpL^0 \to \rpL^1\right)
$$
of $\rpL^0 \cong \pLie$ equals the operad $\Lie$,
\begin{equation}
\label{.,}
H^0(\rpL^*) \cong \Lie.
\end{equation}
The main result of this paper, equivalent to isomorphism~(\ref{.,}), is
therefore a step towards a solution of Problem~\ref{-/}.
}
\end{odstavec}

\begin{odstavec}
\label{gen}
{\rm
{\it Generalizations.\/}
Let us slightly reformulate the above reflections and indicate
possible generalizations. Let $\calP$ be a quadratic 
Koszul operad~\cite[Section~II.3.3]{markl-shnider-stasheff:book} and
$\calB_{\calP_\infty} = (\calB_{\calP_\infty},d)$ the dg-operad of
natural operations on the complex defining the operadic cohomology of
$\calP_\infty$ (= strongly homotopy
$\calP$-algebras~\cite[Definition~II.3.128]{markl-shnider-stasheff:book}) 
with coefficients in itself. In~\cite{markl:de} we conjectured that 
\begin{equation}
\label{oooo}
H^0(\calB^*_{\calP_\infty})
\cong \Lie
\end{equation} 
for each quadratic Koszul operad $\calP$.

The operad $\calB^*_{\calP_\infty}$ has a suboperad
$\calS^*_{\calP_\infty}$ generated by a restricted class of operations
which generalize the braces on the Hochschild cohomology complex of an
associative algebra~\cite{gerstenhaber-voronov:FAP95}. The operad
$\calS^0_{\calP_\infty}$ of degree~$0$ elements in
$\calS^*_{\calP_\infty}$ always contains the operad $\Lie$ for Lie
algebras that represents the intrinsic brackets. 
The conjectural isomorphism~(\ref{oooo}) would therefore imply: 

\begin{conjecture}
\label{jak_to_dopadne}
For each quadratic Koszul operad $\calP$,
$$
\Lie \cong \Ker\left(d : \calS^0_{\calP_\infty} \to \calS^1_{\calP_\infty}
\right) .
$$
\end{conjecture}

Moving from operads to free
algebras~\cite[Section~II.1.4]{markl-shnider-stasheff:book}, an
affirmative solution of this conjecture for a particular operad $\calP$
would immediately give a characterization of Lie elements in free
$\calS^0_{\calP_\infty}$-algebras.

{}From this point of view, the main result of this paper
(Theorem~\ref{.}) is a combination of a solution of
Conjecture~\ref{jak_to_dopadne} for $\calP=\Lie$ with the
identification of $\calS^0_{\Lie_\infty} \cong \pLie$ which expresses
the equivalence between symmetric brace algebras and pre-Lie
algebras~\cite{guin-oudom,lm:sb}.  Conjecture~\ref{jak_to_dopadne}
holds also for $\calP = \Ass$, the operad for associative algebras, as
we know from the Deligne conjecture in the form proved
in~\cite{kontsevich-soibelman}. Since $\calS^0_{\Ass_\infty}$ is the
operad for (ordinary, non-symmetric)
braces~\cite{gerstenhaber-voronov:FAP95}, one can obtain a description of
{\em Lie elements\/} in {\em free brace algebras\/}.  }
\end{odstavec}

\begin{odstavec}
{\rm
{\em Loday's triplettes.\/}
Theorem~\ref{.} can also be viewed as an analog of the characterization of Lie
elements in the tensor algebra $T(V)$ as primitives of the
bialgebra $\calH = (T(V),\otimes,\Delta)$ with $\Delta$ the shuffle
diagonal; we recall this classical result as Theorem~\ref{psano_v_aute} of
Section~2.  The bialgebra $\calH$ is associative,
coassociative cocommutative and its primitives $\Prim(\calH)$ form a
Lie algebra. To formalize such situations, J.-L.~Loday introduced
in~\cite{loday:slides} the notion of a {\em triplette\/}
$(\calC,\spin,\Aalg \stackrel{F}\to \Palg)$, 
abbreviated $(\calC, \calA,\calP)$,
consisting of operads $\calC$ and $\calA$, `spin' relations $\spin$ 
between $\calC$-coalgebras and $\calA$-algebras defining
$(\calC,\spin,\calA)$-bialgebras, an operad $\calP$ describing the
algebraic structure of the primitives, and a forgetful functor $F : \Aalg
\to \Palg$, see
Definition~\ref{jeste_ji_musim_napsat_SMS} in Subsection~\ref{trip}.

The nature of associative, cocommutative coassociative bialgebras and
their primitives is captured by the triplette $(\Com,\Ass,\Lie)$. The
classical Theorem~\ref{psano_v_aute} then follows from the fact that
the triplette $(\Com,\Ass,\Lie)$ is {\em good\/}, in the sense which
we also recall in Subsection~\ref{trip}.  An interesting question is
whether the case of Lie elements in pre-Lie algebras considered in
this paper is governed by a good triplette in which $\calA = \pLie$
and $\calP = \Lie$. See Subsection~\ref{trip} for more detail.  
}
\end{odstavec}

\vskip .3em
\noindent
{\bf Acknowledgments.} I would like to express my thanks to
F.~Chapoton, M.~Livernet, \hbox{J.-L.~Loday}, C.~L\"ofwall,
J.~Stasheff and D.~Tamarkin for many useful comments and
suggestions. I am also indebted to M.~Goze and E.~Remm for their
hospitality during my visit of the University of Mulhouse in the Fall
of 2004 when this work was initiated.

\section{Classical results revisited}
\label{sec1}

In this section we recall some classical results about Lie elements in
free associative algebras in a language suitable for the purposes of
this paper.  Let $\T(V)$ be the tensor algebra generated by a vector
space $V$,
$$
\T(V) = \bfk \oplus \bigoplus_{n=1}^\infty \T^n(V), 
$$
where $\T^n(V)$ is the $n$-th tensor power $\bigotimes^n (V)$ of the
space $V$. Let $\T(V)_{L}$ denote the space $\T(V)$ considered as a Lie
algebra with the commutator bracket
$$
[x,y] := x\ot y - y\ot x,\hskip 1em  x,y \in \T(V),
$$ 
and let $\fL(V) \subset \T(V)$ be the Lie sub-algebra of
$\T(V)_{L}$ generated by $V$. It is well-known that $\fL(V)$ is
(isomorphic to) the
free Lie algebra generated by $V$~\cite[\S4, Theorem~2]{serre:65}. 

There are several characterizations of the subspace $\fL(V)
\subset \T(V)$~\cite{ree:AnM69,serre:65}. 
Let us recall the one which uses the {\em
shuffle diagonal\/} $\Delta : \T(V) \to \T(V) \ot \T(V)$ given, for
$v_1 \ot \cdots \ot v_n \in \T^n(V)$,~by
\begin{equation}
\label{jake_bude_pocasi}
\Delta(v_1 \ot \cdots \ot v_n) := 
\sum_{i=0}^n\sum_{\sigma \in\Sh(i,n-i)} 
[v_{\sigma(1)} \ot \cdots \ot v_{\sigma(i)} ]
\ot [v_{\sigma(i+1)} \ot \cdots \ot v_{\sigma(n)}],
\end{equation}
where $\Sh(i,n-i)$ denotes the set of all $(i,n-i)$-shuffles,
i.e.~permutations $\sigma \in \Sigma_n$ such that
$$
\sigma(1) < \cdots < \sigma(i)\ \mbox { and }  
\sigma(i+1) < \cdots < \sigma(n).
$$

Notice that, in the right hand side of~(\ref{jake_bude_pocasi}), the
symbol $\ot$ has two different meanings, the one inside the brackets
denotes the tensor product in $\T(V)$, the middle one the tensor
product of two copies of $\T(V)$. To avoid this ambiguity, we denote
the product in $\T(V)$ by the dot~$\ccdot$,~(\ref{jake_bude_pocasi})
will then read as
$$
\Delta(v_1 \ccdot \cdots \ccdot v_n) := 
\sum_{i=0}^n\sum_{\sigma \in \Sh(i,n-i)} 
[v_{\sigma(1)} \ccdot \cdots \ccdot v_{\sigma(i)}] 
\ot [v_{\sigma(i+1)} \ccdot\cdots \ccdot v_{\sigma(n)}].
$$

The triple $(\T(V),\ccdot,\Delta)$ is a standard example of a unital counital
associative coassociative cocommutative Hopf algebra. 
We will need also the {\em augmentation ideal\/}
$\bfA(V) \subset \fA(V)$ which equals $\fA(V)$ minus the ground field, 
$$
\bfA(V) = \bigoplus_{n=1}^\infty \T^n(V),
$$
and the {\em reduced\/} diagonal $\bDelta : \bfA(V) \to \bfA(V) \ot
\bfA(V)$ defined as 
$$
\bDelta(x) := \Delta(x) - 1 \ot x - x \ot 1,\ \mbox {  for $x \in \bfA(V)$,}
$$
or, more explicitly,
$$
\bDelta(v_1 \ccdot \cdots \ccdot v_n) := 
\sum_{i=1}^{n-1}\sum_{\sigma \in \Sh(i,n-i)} 
[v_{\sigma(1)} \ccdot \cdots \ccdot v_{\sigma(i)}] 
\ot [v_{\sigma(i+1)} \ccdot\cdots \ccdot v_{\sigma(n)}],
$$
for $v_1,\ldots,v_n \in V$ and $n \geq 1$. Clearly $\L(V) \subset
\bT(V)$. The following theorem is classical~\cite{serre:65}.

\begin{theorem}
\label{psano_v_aute}
The subspace $\fL(V) \subset \bT(V)$ equals the subspace of
primitive elements,
$$
\fL(V) = \Ker\left(\bDelta : \bfA(V) \to \bfA(V)\ot \bfA(V)\right).
$$
\end{theorem}

The diagonal $\Delta : \T(V) \to \T(V) \ot \T(V)$ is a homomorphism of
associative algebras, that~is 
\begin{equation}
\label{to_jsem_zvedav_jesli_budou_pastviny_fungovat}
\Delta (x \ccdot y) = \Delta(x) \ccdot \Delta(y),\ \mbox {  for $x \in \T(V)$,}
\end{equation}
where the same $\ccdot$ denotes both the multiplication in $\T(V)$ in
the left hand side and the induced multiplication of $\T(V) \ot \T(V)$
in the right hand side. The reduced diagonal $\bDelta : \bfA(V) \to
\bfA(V) \ot \bfA(V)$ is, however, of a different nature:

\begin{proposition}
For each $x,y \in \bfA(V)$,
\begin{equation}
\label{1}
\bDelta(x\ccdot y) = \Delta(x)\ccdot\bDelta(y) + \bDelta(x)\ccdot(y \ot 1 +
1 \ot y) + (x \ot y + y \ot x).
\end{equation}
\end{proposition}

The proof is a direct verification which we leave for the reader.  We
are going to reformulate~(\ref{1}) using an action of $\bT(V)$ on
$\bT(V) \ot \bT(V)$ defined as follows. For $\xi \in \bfA(V) \ot
\bfA(V)$ and $x \in \bfA(V)$, let
\begin{equation}
\label{2}
\begin{array}{rcl}
x* \xi &:=& \Delta(x)\ccdot  \xi \in \bfA(V) \ot \bT(V),\ \mbox { and}
\\
\xi* x &:=& \xi \ccdot  (1 \ot x + x \ot 1)  \in \bfA(V)  \ot \bT(V),
\rule{0em}{1.2em}
\end{array}
\end{equation}
where $\ccdot$ denotes, as before, the tensor multiplication in
$\bT(V) \ot \bT(V)$.  Observe that, while
\begin{equation}
\label{3}
(x\ccdot y)* \xi = x *(y* \xi)\ \mbox { and }\ (x*\xi)*y = x*(\xi* y),
\end{equation}
$(\xi* x)* y \not=\xi* (x\ccdot y)$, therefore the action~(\ref{2})
{\em does not\/} make $\bT(V) \ot \bT(V)$ a bimodule over the associative
algebra $(\bfA(V),\ccdot)$.  To understand the algebraic properties of
the above action better, we need to recall the following important

\begin{Definition}{\rm (\cite{gerstenhaber:AM63})}
A {\em pre-Lie\/} algebra is a vector space $X$ with a bilinear
product $\star : X \ot X \to X$ such that the associator $\Phi :
X^{\ot 3} \to X$ defined by
\begin{equation}
\label{dnes_jsem_objel_Polednik}
\Phi(x,y,z):=(x \star y) \star z-x\star (y\star z),\ \mbox { for
$x,y,z \in X$},
\end{equation}
is symmetric in the last two variables,
$\Phi(x,y,z) = \Phi(x,z,y)$. Explicitly,
\begin{equation}
\label{aby_to_nezkoncilo_prusvihem}
(x\star y) \star z-x\star (y\star z) = 
(x\star z) \star y-x\star (z\star y)\ \mbox { for each $x,y,z \in X$}.
\end{equation}
\end{Definition}

There is an obvious graded version of this definition.  Pre-Lie
algebras are known also under different names, such as right-symmetric
algebras, Vinberg algebras, \&c.  Pre-Lie algebras are particular
examples of {\em Lie-admissible\/} algebras~\cite{markl-remm}, which means
that the object $X_L := (X,[-,-])$ with $[-,-]$ the commutator of
$\star$, is a Lie algebra.  Each associative algebra is clearly
pre-Lie. In the following proposition, $\bfA(V)_{pL}$ denotes the
augmentation ideal $\bT(V)$ of the associative algebra $\T(V)$
considered as a pre-Lie algebra.

\begin{proposition}
Formulas~(\ref{2}) define on $\bfA(V) \ot \bfA(V)$ a structure of a
bimodule over the pre-Lie algebra $\bfA(V)_{pL}$. This means that
$$
(\xi* x)*y - \xi * (x \ccdot y) = (\xi * y) * x - \xi* (y \ccdot x)
$$
and
$$
(x \ccdot y) * \xi - x * (y*  \xi) = (x * \xi) *y - x*(\xi* y),
$$
for each $x,y \in \bT(V)$ and $\xi \in \bT(V) \ot \bT(V)$. In
particular, $\bfA(V) \ot \bfA(V)$ is a module over the Lie algebra
$\bfA(V)_L$.
\end{proposition}

\begin{Proof}
To prove the first equality, notice that
$$
(\xi * x) *y - \xi *(x \ccdot y) 
= \xi\ccdot(x \ot y + y \ot x) =  (\xi* y)*x - \xi*(y \ccdot x).
$$
The second one immediately follows from~(\ref{3}).
\end{Proof}

Using action~(\ref{2}), rule~(\ref{1}) can be rewritten as
\begin{equation}
\label{4}
\bDelta(x \ccdot y) = \bDelta(x)*y + x*\bDelta(y) +  R(x,y),\
x,y \in \bT(V),
\end{equation}
where the symmetric bilinear form $R(x,y) := x \ot y + y \ot x$
measures the deviation of $\bDelta$ from being a pre-Lie algebra
derivation in 
$$
{\it Der}_{\it pre-Lie}\left(\bfA(V)_{pL},\bfA(V) \ot
\bfA(V)\right).
$$

On the other hand, since $R: \bT(V) \to \bT(V) \ot \bT(V)$ is
symmetric, $\bDelta$ is a derivation of the associated Lie algebra
$\bfA(V)_L$,
$$
\bDelta \in {\it Der}_{\it Lie}\left(\bfA(V)_{L},\bfA(V) \ot \bfA(V)\right),
$$
which implies that $\L(V) \subset \Ker(\bDelta)$. The following
statement is completely obvious and we formulate it only to motivate
Proposition~\ref{a} of Section~3.

\begin{proposition}
\label{strasne_pocasi_trva}
The map $\bDelta : \bT(V) \to \bT(V) \ot \bT(V)$ is uniquely
determined by the rule~(\ref{4}) together with the requirement that
$\bDelta(v) = 0$ for $v \in V$.
\end{proposition}

Observe that the reduced diagonal $\bDelta : \bfA(V) \to \bfA(V) \ot
\bfA(V)$ is the initial differential of the cobar construction
\begin{equation}
\label{1a}
{\it Cob}(\bT(V),\bDelta): \hskip 1em
\bfA(V) \stackrel{d}{\longrightarrow}
\bfA(V) \ot \bfA(V)  \stackrel{d}{\longrightarrow}
\bfA(V)\ot \bfA(V)\ot \bfA(V) \stackrel{d}{\longrightarrow}\cdots
\end{equation}
of the coassociative coalgebra $(\bfA(V),\bDelta)$.
Complex~(\ref{1a}) calculates the cohomology
\begin{equation}
\label{ooOoo}
H^*\left(\T(V),\Delta\right) = 
{\it Cotor}^{*+1}_{(\T(V),\Delta)}(\bfk,\bfk)
\end{equation} 
of the shuffle coalgebra.  

On the other hand, 
by the Poincar\'e-Birkhoff-Witt theorem,
there is an isomorphism of coalgebras
$$
(\T(V),\Delta) \cong (\bfk[\L(V)],\nabla),
$$
where the polynomial ring $\bfk[\L(V)]$ in the right hand side is equipped
with the standard cocommutative comultiplication $\nabla$. 
Dualizing the proof of the
classical~\cite[Theorem~VII.2.2]{maclane:homology},
one obtains the isomorphism
$$
{\it Cotor}^{*}_{(\bfk[\L(V)],\nabla)}(\bfk,\bfk) \cong \ext^*(\L(V))
$$ 
where $\ext^*(-)$
denotes the exterior algebra functor.
We conclude that 
$$
H^*\left(\T(V),\Delta\right) \cong \ext^{*+1}(\L(V)).
$$

\section{Lie elements in the free pre-Lie algebra}
\label{sec2}

In this section we show that the results reviewed in
Section~2 translate to pre-Lie algebras.
Let $\pL(V) = (\pL(V),\star)$ denote the free pre-Lie algebra
generated by a vector space $V$ and let $\pL(V)_L$ be the associated
Lie algebra. The following proposition is proved
in~\cite{dzhu-lof:HHA02}, but we will give a shorter and more
direct proof, which was kindly suggested to us by M.~Livernet.

\begin{proposition}
The subspace $\L(V) \subset \pL(V)_L$ generated by $V$ is isomorphic
to the free Lie algebra on $V$.
\end{proposition}

\noindent
{\bf Proof}
(due to M.~Livernet).\hglue 1.8em 
Let us denote in 
this proof by $\L'(V)$ the Lie subalgebra of $\pL(V)_L$ generated by $V
\subset \pL(V)$ and by $\L''(V)$ the Lie subalgebra of $\T(V)_L$ generated by
$V \subset \T(V)$. The canonical map $\pL(V) \to \T(V)_{pL}$ clearly
induces a map $\pL(V)_L \to \T(V)_L$ which restricts to a Lie algebra
homomorphism $\alpha: \L'(V) \to \L''(V)$.

Let $\L(V)$ be, as before, the free Lie algebra generated by
$V$. Since $\L'(V)$ is also generated by $V$, the canonical map $\beta
: \L(V) \to \L'(V)$ is an epimorphism. To prove that it is a
monomorphism, observe that the composition $\alpha \beta : \L(V) \to
\L''(V)$ coincides with the canonical map induced by the inclusion $V
\hookrightarrow \L''(V)$. Since $\L''(V)$ is isomorphic to the free
Lie algebra generated by $V$~\cite[\S4, Theorem~2]{serre:65}, the
composition $\alpha \beta : \L(V) \to \L''(V)$ is an isomorphism,
therefore $\beta$ must be monic. We conclude that the canonical map
$\beta : \L(V) \to \L'(V)$ is an isomorphism, which finishes the
proof.\qed

Consider the free graded pre-Lie algebra $\pL(V,\oh)$ generated by $V$
and one `dummy' variable $\oh$ placed in degree $+1$.  Observe that
\begin{equation}
\label{5}
\pL^*(V,\oh) = \pL(V) \oplus \bigoplus_{n \geq 1} \pL^n(V), 
\end{equation}
where $\pL^n(V)$ is the subset of $\pL(V)$ spanned by
monomials with exactly $n$ occurrences of the dummy variable $\oh$.

We need to consider also the graded pre-Lie algebra $\rpL(V)$ (``r''
for ``reduced'') defined as the quotient
$$
\rpL(V) := \pL(V,\oh)/(\oh \star \oh)
$$
of the free pre-Lie algebra $\pL(V,\oh)$ by the ideal $(\oh \star \oh)$
generated by $\oh \star \oh$.  The grading~(\ref{5}) clearly induces a
grading of $\rpL(V)$ such that $\rpL^0(V) = \pL(V)$ and $\rpL^1(V) =
\pL^1(V)$,
\begin{equation}
\label{koupali-jsme-se-na-Hradistku}
\rpL^*(V) = \pL(V) \oplus \pL^1(V) \oplus \bigoplus_{n \geq 2}
\rpL^n(V).
\end{equation}
The following statement, in which $\Phi$ is the
associator~(\ref{dnes_jsem_objel_Polednik}), is an analog of
Proposition~\ref{strasne_pocasi_trva}.

\begin{proposition}
\label{a}
There exists precisely one degree $+1$ map $d : \rpL^*(V) \to
\rpL^{*+1}(V)$ such that $d(v) = 0$ for $v \in V$, $d(\oh) = 0$ and
\begin{equation}
\label{0}
d(a \star b) = d(a) \star b + (-1)^{|a|} a \star d(b) + Q(a,b),
\end{equation}
where 
\begin{equation}
\label{6}
Q(a,b) : = (\oh \star a ) \star b - \oh \star (a \star b) = \Phi(\oh,a,b),
\end{equation}
for $a,b \in \rpL^*(V)$. Moreover, $d^2 = 0$.
\end{proposition}

\begin{Proof}
The uniqueness of the map $d$ with the properties stated
in the proposition is clear.  To prove that such a map exists, we show
first that there exists a degree one map $\td : \pL(V,\oh) \to
\pL(V,\oh)$ of graded free pre-Lie algebras such that $\td(v) = 0$ for $v
\in V$, $\td(\oh) = 0$ and
\begin{equation} 
\label{0'}
\td(x \star y) = \td(x) \star y + (-1)^{|x|} x \star \td(y) + Q(x,y),
\end{equation}
where $Q(x,y) : = \Phi(\oh,x,y)$ for $x,y \in \pL(V,\oh)$. Let us
verify that the above rule is compatible with the axiom $\Phi(x,y,z) =
(-1)^{|z||y|}\Phi(x,z,y)$ of graded pre-Lie
algebras. Applying~(\ref{0'}) twice, we obtain
\begin{eqnarray}
\label{16a}
\td\Phi(x,y,z) \!\! &=&\!\!
\Phi(\td x,y,z) + (-1)^{|x|} \Phi(x,\td y,z) + (-1)^{|x| +
  |y|}\Phi(x,y,\td z) 
\\ \nonumber 
&&
- (-1)^{|x|} x \star Q(y,z) + Q(x \star y,z) +Q(x,y) \star z -
Q(x,y \star z),
\end{eqnarray}
for arbitrary $x,y,z \in \pL(V,\oh)$.  

Let us make a small digression and observe that the associator $\Phi$
behaves as a Hochschild cochain, that is
$$
\oh \star\Phi (x,y,z) - \Phi(\oh \star x,y,z) + \Phi(\oh, x \star y,z)
- \Phi(\oh,x,y\star z) + \Phi(\oh,x,y) \star z = 0.
$$
It follows from the definition of the form $Q$ and the above equation
that the last three terms of~(\ref{16a}) equal $\Phi(\oh \star x,y,z)
- \oh \star\Phi (x,y,z)$, therefore~(\ref{16a}) can be rewritten as 
\begin{eqnarray*}
\td\Phi(x,y,z) &=&
\Phi(\td x,y,z) + (-1)^{|x|} \Phi(x,\td y,z) + (-1)^{|x| +
  |y|}\Phi(x,y,\td z) 
\\ \nonumber 
&&
- (-1)^{|x|} x \star Q(y,z) + \Phi(\oh \star x,y,z)
- \oh \star\Phi (x,y,z).
\end{eqnarray*}
Since the right hand side of the
above equality is graded symmetric in $y$ and $z$, we conclude that
$$
\td \left(\Phi(x,y,z) - (-1)^{|z||y|}\Phi(x,z,y)\right) =0,
$$
which implies the existence of $\td : \pL(V,\oh) \to \pL(V,\oh)$ with
the properties stated above.  It is easy to verify, using~(\ref{0'})
and the assumption $\td(\oh) = 0$, that
\begin{equation}
\label{k}
\td(\oh \star \oh) = \Phi(\oh,\oh,\oh)
\end{equation}
and that
\begin{equation}
\label{l}
\td^2(x \star y) = \td^2(x) \star y + x \star \td^2(y) + Q(\td x,y) +
(-1)^{|x|} Q(x,\td y) + \Phi(\oh \star \oh,x,y)
\end{equation}
for arbitrary $x,y \in \pL(V,\oh)$. 

A simple induction on the number of generators based on~(\ref{k})
together with the rule~(\ref{0'}) shows that $\td$ preserves the ideal
generated by $\oh \star \oh$. An equally simple induction based
on~(\ref{l}) and~(\ref{0'}) shows that $\Im(\td^2)$ is a subspace of
the same ideal. We easily conclude from the above facts that $\td$
induces a map $d : \rpL^*(V) \to \rpL^{*+1}(V)$ required by the
proposition.
\end{Proof}

Let us remark that each pre-Lie algebra
$(X,\star)$ determines a unique {\em symmetric brace algebra\/}
$(X,-\langle -,\ldots,-\rangle)$ with $x \langle y \rangle = x\star y$
for $x,y \in X$~\cite{lm:sb,guin-oudom}. 
The bilinear form $Q$ in~(\ref{0}) then can be
written as
$$
Q(a,b) = \oh \langle a,b \rangle,\ \mbox { for } a,b \in \rpL(V).
$$
The complex
\begin{equation}
\label{za_chvili_tam_musim_volat_a_hrozne_se_mi_nechce!}
\pL(V) \stackrel{d}{\longrightarrow} \pL^1(V)
\stackrel{d}{\longrightarrow} \rpL^2(V) \stackrel{d}{\longrightarrow} \cdots 
\end{equation}
should be viewed as an analog of the cobar construction~(\ref{1a}).
We will see in Section~5 that it describes natural
operations on the Chevalley-Eilenberg cohomology of a Lie algebra. 
The main result of this paper reads:

\begin{theorem}
\label{.}
The subspace $\fL(V) \subset \pL(V)$ equals the kernel of the map $d :
\pL(V) \to \pL^1(V)$,
$$
\fL(V) = \Ker \left(\rule{0em}{1em} d : \pL(V) \to \pL^1(V)  \right).
$$
\end{theorem}

In Section~4 we
describe the spaces $\pL(V)$, $\pL^1(V)$ and the map $d : \pL(V) \to
\pL^1(V)$ in terms of trees.  Theorem~\ref{.} will be proved in
Section~6.

\section{Trees}
\label{v_podstate_cela_napsana_v_Moravske_Trebove}

We begin by recalling a tree description of free pre-Lie algebras due
to F.~Chapoton and M.~Livernet~\cite{chapoton-livernet:pre-lie}. By a
{\em tree\/} we understand a finite connected simply connected graph
without loops and multiple edges. We will always assume that our trees
are {\em rooted\/} which, by definition, means that one of the
vertices, called the {\em root\/}, is marked and all edges are
oriented, pointing to the root.

Let us denote by $\Tr_n$ the set of all trees with $n$ vertices
numbered $1,\ldots,n$. The symmetric group $\Sigma_n$ act on
$\Tr_n$ by relabeling the vertices. We define
$$
\Tr_n(V) := {\rm Span}_{\bfk}{(\Tr_n)} \ot_{\Sigma_n} \otexp Vn,\
n \geq 1,
$$
where ${\rm Span}_{\bfk}{(\Tr_n)}$ denotes the $\bfk$-vector space
spanned by $\Tr_n$ with the induced $\Sigma_n$-action and where
$\Sigma_n$ acts on $\otexp Vn$ by permuting the
factors. Therefore $\Tr_n(V)$ is the set of trees with $n$
vertices decorated by elements of $V$.

\begin{example}
The set $\Tr_1$ consists of a single tree 
\hskip .3em $\bullet$ \hskip .3em with one vertex (which is also the
root) and no edges, thus $\Tr_1(V) \cong V$.  The set $\Tr_2$
consists of labelled trees
$$
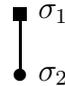
\stromI{\sigma_1}{\sigma_2}
$$
where $\hskip .1em \root \hskip .1em $ 
denotes the root and $\sigma \in \Sigma_2$. This means
that $V$-decorated trees from $\Tr_2(V)$ look~as
$$
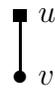
\stromI uv
$$
where $u,v \in V$, therefore $\Tr_2(V) \cong \otexp V2$.  Similarly,
$\Tr_3(V) \cong \otexp V3 \oplus (V \ot S^2(V))$, where $S^2(V)$
denotes the second symmetric power of $V$. The corresponding decorated
trees are
$$
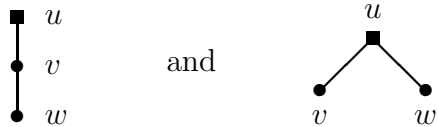
\rule{0em}{4em}
\mbox{\stromII uvw \hskip 3em  \raisebox{1.5em}{ and }\hskip 3em \stromIII uvw}
$$
for $u,v,w \in V$. Finally,
$$
\Tr_4(V) \cong \otexp V4 \oplus \otexp V4 \oplus (\otexp V2 \ot S^2(V))
\oplus (V \ot S^3(V)),
$$
with the summands corresponding to the decorated trees
\begin{equation}
\label{snad_nenachladnu}
\raisebox{-2.7em}{\rule{0em}{5em}}
\stromIV tuvw \hskip 3em   \stromV tuvw  \hskip 3em   \stromVI tuvw 
 \hskip 3em   \stromVII tuvw 
\end{equation}
with $t,u,v,w \in V$.
\end{example}

\begin{theorem}
{\rm (Chapoton-Livernet~\protect\cite{chapoton-livernet:pre-lie})}
\label{jsem_zvedav_jak_povezu_tu_VOSu_v_desti}
Let $\Tr(V) := \bigoplus_{n \geq 1}\Tr_n(V).$ Then there is a natural
isomorphism
\begin{equation}
\label{zitra_odjizdim_na_zavody}
\pL(V) \cong \Tr(V).
\end{equation}
\end{theorem}
 
The pre-Lie multiplication in the left hand side
of~(\ref{zitra_odjizdim_na_zavody}) translates to the vertex
insertion of decorated trees in the right hand side,
see~\cite{chapoton-livernet:pre-lie} for details.

\begin{example}
The most efficient way to identify decorated trees with elements
of free pre-Lie algebras is to use the formalism of symmetric brace
algebras~\cite{lm:sb}.  The trees in~(\ref{snad_nenachladnu})
then represent the following elements of $\pL(V)$:
$$
t \langle u  \langle v  \langle w  \rangle \rangle \rangle,\
t  \langle u,v  \langle w  \rangle \rangle, \
t  \langle u  \langle v,w  \rangle  \rangle\
\mbox { and } \ t  \langle u,v,w \rangle.
$$
\end{example}

Using the same tree description~\cite{chapoton-livernet:pre-lie} of the
free graded pre-Lie algebra $\pL(V,\oh)$, one can easily get
a natural isomorphism
\begin{equation}
\label{zitra_odjizdim_na_zavody-uz_tam_sem}
\pL^1(V) \cong \Tr^1(V) :=  \bigoplus_{n \geq 0} \Tr_n^1(V),
\end{equation}
where $\Tr_n^1(V)$ is the set of all trees with $n$ vertices decorated
by elements of $V$ and one vertex decorated by the dummy variable
$\oh$. We call the vertex decorated by $\oh$ the {\em special
vertex\/}.

\begin{example}
Clearly $\Tr_0^1(V) \cong \bfk$ while $\Tr_1^1(V) \cong V \oplus V$ with
the corresponding decorated trees
$$
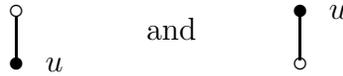
\stromVIII \oh\bullet {}u 
\hskip 3em  \raisebox{.8em}{ and }\hskip 3em 
\stromVIII \bullet\oh u{} 
$$
where $u \in V$. Similarly, 
$$
\Tr_2^1(V) \cong \otexp V2 \oplus \otexp V2 \oplus \otexp V2 \oplus 
\otexp V2 \oplus  S^2(V)
$$
with the corresponding trees
$$
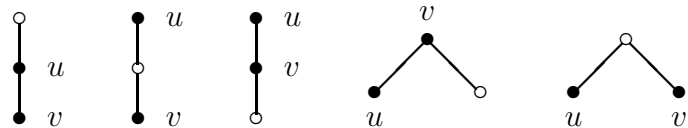
\stromIX \oh\bullet\bullet{}uv
\hskip 3em 
\stromIX \bullet\oh\bullet u{}v
\hskip 3em 
\stromIX \bullet\bullet\oh uv{}
\hskip 3em 
\stromX \bullet\bullet\oh uv{}
\hskip 3em 
\stromX \bullet\oh\bullet u{}v
$$
for $u,v \in V$.
In the above pictures we always placed the root
on the top. Some examples of decorated trees from
$\Tr_n^1(V)$, $n\geq 3$, can also be found in
Examples~\ref{abych_jeste_neonemocnel},
\ref{uz_ctvrty_neletovy_den_v_rade}, \ref{mmm} and~\ref{yyy}.
\end{example}

Let us describe the map $d : \pL(V) \to \pL^1(V)$ of Theorem~\ref{.}
in terms of decorated trees.  We say that an edge $e$ of a decorated
tree $S\in \Tr_n^1(V)$ is {\em special\/} if it is adjacent to the
special vertex of $S$. Given such an edge $e$, we define
the quotient $S/e \in \Tr_n(V)$ by contracting the special edge of $S$
into a vertex and decorating this vertex by the label of the (unique)
endpoint of $e$ different from the special vertex.  In the following
examples, the special edge will be marked by the double line.

\begin{example}
\label{abych_jeste_neonemocnel}
If $S \in \Tr^1_3(V)$ is the tree
$$
\raisebox{.5em}{S \hskip .3em = \hskip .3em \stromXIVe  \hskip .3em ,}
$$
$u,v,w \in V$, then
\vskip .3em 
$$
S/e \hskip .3em = \hskip .6em  \raisebox{-1em}{\stromIII vuw}\hskip .3em .
$$
\end{example}

Let $T \in \Tr_n(V)$. We call a couple $(S,e)$, where $S \in
\Tr_n^1(V)$ and $e$ a special edge of $S$, a {\em blow-up\/} of\/
$T$ if $S/e \cong T$ and if the arity (= the number of incoming
edges) of the special vertex of $S$ is $\geq 2$.
We denote by $\Bl(T)$ the set of all blow-ups of $T$.

\begin{example}
\label{uz_ctvrty_neletovy_den_v_rade}
The set $\Bl(\hskip .1em \raisebox{.1em}{\root}\hskip .1em)$ is
empty. The simplest nontrivial example of a blow-up is
\begin{eqnarray*}
\Bl \left(\hskip -2.5em  
\raisebox {-.8em}{\stromI uv} \hskip .3em\right) &=&   
\left\{ 
\stromXI \bullet\oh\bullet u{}v 
\right\}, 
\end{eqnarray*}
\vskip .2em
\noindent 
where the double line denotes, as in
Example~\ref{abych_jeste_neonemocnel}, the special edge. Let us give
two more examples where $u$, $v$ and $w$ are elements of $V$:
\begin{eqnarray*}
\Bl \left(\rule{0pt}{2.7em}\hskip 1em  \raisebox {-1.3em}{\stromII uvw} \hskip
.6em\right) &=&
\left\{ \rule{0em}{3em} \hskip .4em
\stromXII \hskip .4em ,  \stromXIII\hskip .4em
\right\} \hskip .4em \mbox { and }
\\
\Bl \left(\hskip .4em  \raisebox {-1.3em}{\stromIII uvw} \hskip
.4em\right) &=&
\left\{ \rule{0em}{3em} \hskip .4em
\stromXIV \hskip .4em , \hskip .4em
\stromXV \hskip .4em ,  \stromXVI\hskip .4em, \stromeXtra\hskip .4em
\right\}.
\end{eqnarray*}
\end{example}

The last thing we need is to introduce, for $(S,e) \in \Bl(T)$, the 
sign $\epsilon_{(S,e)} \in \{-1,+1\}$ as
$$
\epsilon_{(S,e)} := \cases {+1}{if $e$ is an
 incoming edge of the special vertex, and}
{-1}{if $e$ is the
outgoing edge of the special vertex.}
$$
Finally, define the map
\begin{equation}
\label{Nunyk}
\delta : \Tr(V) \to \Tr^1(V)
\end{equation}
by
$$
\delta(T) : = \sum_{(S,e) \in \Bl(T)} \epsilon_{(S,e)} S.
$$

\begin{example}
\label{mmm}
In this example, $t,u,v$ and $w$ are arbitrary elements of $V$. We
stick to our convention that the root is placed on the top. 
Let us give first some examples of the map $\delta : \Tr(V) \to
\Tr^1(V)$ that follow immediately from the calculations in
Example~\ref{uz_ctvrty_neletovy_den_v_rade}. We keep the 
double lines indicating
which edges has been blown-up:
\begin{eqnarray*}
\delta(\hskip .1em \raisebox{.1em}{\root}\hskip .1em ) &=& 0,
\\
\delta\left(\hskip -2.5em  
\raisebox {-.8em}{\stromI uv} \hskip .4em
\right) 
&=& \stromXI \bullet\oh\bullet u{}v, \rule{0em}{2.5em}
\\
\rule{0em}{3.3em}\delta\left( \hskip -.4em \rule{0em}{2.3em} \hskip 1em  
\raisebox{-1.3em}{\stromII uvw}\hskip .8em  \right) &=&  
\stromXII \hskip .4em + \stromXIII,
\\
\rule{0em}{3.3em}\delta \left(\hskip .4em  
\raisebox {-1.3em}{\stromIII uvw} \hskip .4em\right) 
&=& 
\stromXV  \hskip .4em +  \hskip .4em \stromXIV  \hskip .4em
+  \hskip .4em \stromeXtra  \hskip .4em -
\stromXVI.
\end{eqnarray*}

\vskip .5em
\noindent 
Let us give some more formulas, this time without indicating the
blown-up edges:
\begin{eqnarray*}
\rule{0em}{3.7em}\delta 
\left(\hskip 1.3em\raisebox{.4em}{\stromIV tuvw} \right) &=&
\stromXX + \stromXXI + \stromXXII,
\\
\rule{0em}{3.7em}\delta 
\left(\hskip .3em\raisebox{.4em}{\stromV tuvw} \hskip .1em \right) &=&
\stromXXIII + \stromXXIV + \stromXXV 
\\
&& \hskip 2em  \raisebox{1.5em}{ + \stromXXVIextra - \stromXXVI,}
\\
\rule{0em}{3.3em}\delta \left(\hskip .4em  
\raisebox {.3em}{\stromVI tuvw} \hskip .4em\right) &=&
\stromXXVII + \stromXXVIII - \stromXXIX
\\
&&  \hskip 2em  \raisebox{1.5em}{+ \stromXXX +\stromExtra.}
\end{eqnarray*}
\end{example}

The proof of the following proposition is a direct verification based
on the induction on the number of vertices and formula~(\ref{0}).

\begin{proposition}
\label{abych_neonemocnel}
The diagram
$$
\unitlength 1em
\thicklines
\begin{picture}(10,7)(7,10)
\put(6,16){\makebox(0,0)[cc]{$\pL(V)$}}
\put(16,16){\makebox(0,0)[cc]{$\pL^1(V)$}}
\put(6,10){\makebox(0,0)[cc]{$\Tr(V)$}}
\put(16,10){\makebox(0,0)[cc]{$\Tr^1(V)$}}
\put(8,16){\vector(1,0){6}}
\put(8,10){\vector(1,0){6}}
\put(6,15){\vector(0,-1){4}}
\put(16,15){\vector(0,-1){4}}
\put(11,17){\makebox(0,0)[cc]{$d$}}
\put(11,11){\makebox(0,0)[cc]{$\delta$}}
\put(7,13){\makebox(0,0)[cc]{$\cong$}}
\put(17,13){\makebox(0,0)[cc]{$\cong$}}
\end{picture}
$$
in which the vertical maps are
isomorphism~(\ref{zitra_odjizdim_na_zavody})
and~(\ref{zitra_odjizdim_na_zavody-uz_tam_sem}), is commutative.
\end{proposition}

\begin{corollary}
\label{,}
There is a natural isomorphism
$$
\L(V) \cong \Ker ( \delta : \Tr(V) \to \Tr^1(V)).
$$
\end{corollary}

\begin{example}
\label{yyy}
It follows from the formulas given in Example~\ref{mmm} that, for each
$u,v,w \in V$,
\def\po#1{{\hskip 1em \raisebox{1.5em}{#1} \hskip .7em }}
\def\poo#1{{\hskip .7em \raisebox{1.5em}{#1} \hskip .5em }}
\begin{eqnarray*}
\lefteqn{
\rule{0em}{3.3em}\delta\left( \hskip -.4em \rule{0em}{2.3em} \hskip 1em  
\raisebox{-1.3em}{
\stromII uvw \po- \stromII uwv \po- \stromII vwu \po+ \stromII wvu
}\hskip .8em  \right) =}
\\
&=& \rule{0em}{2.5em} \hskip .5em  
\stromXii uvw \hskip .5em-\hskip .5em \stromXii uwv\hskip .5em 
-\hskip .5em 
\stromXii vwu \hskip .5em+ \hskip .5em\stromXii  wvu
\\\rule{0em}{4em}
&&+ \stromXvi uvw - \stromXvi uwv - \stromXvi  vwu + \stromXvi wvu 
\\
&=& 
\rule{0em}{5em}\delta \left(\hskip .4em  
\raisebox {-1.3em}{\stromIII vuw\poo- \stromIII wuv} \hskip .4em\right),
\end{eqnarray*}
therefore the combination
$$
\xi_{u,v,w} := 
\raisebox{-1.3em}{
\stromII uvw \po- \stromII uwv \po- \stromII vwu \po+ \stromII wvu
\po- \stromIII vuw \poo+ \stromIII wuv}
$$
belongs to the kernel of $\delta : \Tr_3(V) \to \Tr_3^1(V)$. It is
easy to see that elements of this form in fact span this kernel and
that the correspondence $\xi_{u,v,w} \mapsto [u,[v,w]]$ defines an
isomorphism
$$
\Ker \left(\Tr_3(V) \to \Tr_3^1(V)\right) \cong \L_3(V),
$$
where $\L_3(V) \subset \L(V)$ denotes the subspace of
elements of monomial length~$3$.
\end{example}

\section{Cohomology operations}
\label{tt}

In this section we show how an object closely related to the cochain
complex $\rpL^*(V) = (\rpL^*(V),d)$
of~(\ref{koupali-jsme-se-na-Hradistku}), considered in
Proposition~\ref{a}, naturally acts on the Chevalley-Eilenberg complex
of a Lie algebra with coefficients in itself.  For $n \geq 1$, let
$\bfk^n := \Span_{\bfk}(\Rada e1n)$ and let $\rpL^*(n)$ denote the subspace
of the graded vector space $\rpL^*(\bfk^n)$ spanned by monomials which
contain each basic element $\Rada e1n$ exactly once.

More formally, given an $n$-tuple $\Rada t1n \in \bfk$,
consider the map $\varphi_{\Rada t1n} : \bfk^n \to \bfk^n$ defined~by
$$
\varphi_{\Rada t1n}(e_i) := t_i e_i,\ 1 \leq i \leq n.
$$
Let us denote by the same symbol also the induced map $\varphi_{\Rada
t1n} : \rpL^*(\bfk^n) \to \rpL^*(\bfk^n)$. Then
$$
\rpL^*(n) := \left\{\rule{0pt}{.9em} x \in \rpL^*(\bfk^n);\
\varphi_{\Rada t1n}(x) = t_1\cdots t_n x \ \mbox { for each $\Rada t1n
  \in \bfk$}\right\}. 
$$
The above description immediately implies that $\rpL^*(n)$ is a $d$-stable
subspace of $\rpL^*(\bfk^n)$, therefore $\rpL^*(n) = (\rpL^*(n),d)$ is
a chain complex for each $n\geq 1$.  Clearly $\pL(n) \cong
\Span_{\bfk}(\Tr_n)$ and $\pL^1(n) \cong \Span_{\bfk}(\Tr_n^1)$.
Observe that the above reduction does not erase any information,
because $\rpL^*(V)$ can be reconstructed as
$$
\rpL^*(V) \cong \bigoplus_{n \geq 1}\rpL^*(n) \ot_{\Sigma_n} \otexp Vn.
$$

Let us explain how each $U \in \rpL^d(n)$ determines an
$n$-multilinear degree $d$ operation on the Chevalley-Eilenberg
complex $\CE*$ of a Lie algebra $L$ with coefficients in
itself~\cite{chevalley-eilenberg}. 
We will use the standard
identification~\cite[Definition~II.3.99]{markl-shnider-stasheff:book}
\begin{equation}
\label{potim_se}
\CE * \cong \Coder^*(\L^c(\desusp L))
\end{equation}
where $\L^c(\desusp L)$ denotes the cofree conilpotent Lie
coalgebra~\cite{markl-shnider-stasheff:book} 
cogenerated by the desuspension $\desusp L$ of
the vector space $L$. Let $\lambda \in \Coder^1(\L^c(\desusp L))$
be the co-extension of the desuspended Lie algebra bracket 
$$
\desusp \circ [-,-]\circ (\uparrow \land \uparrow) : \desusp L \land \desusp
L \to \desusp L
$$ 
into a coderivation. Then $\lambda^2=0$ and~(\ref{potim_se})
translates the Chevalley-Eilenberg differential $\dce$ into the
commutator with $\lambda$. 

The above construction can be easily generalized to the case when $L$
is an $L_\infty$-algebra, $L = (L,l_1,l_2,l_3,\ldots)$~\cite{lada-markl:CommAlg95}. The
structure operations $(l_1,l_2,l_3,\ldots)$  assemble again into a
coderivation $\lambda \in \Coder^1(\L^c(\desusp L))$ with
$\lambda^2=0$~\cite[Theorem~2.3]{lada-markl:CommAlg95}, and~(\ref{potim_se}) can be taken for a
definition of the (Chevalley-Eilenberg) cohomology of 
$L_\infty$-algebras with coefficients in itself.

The last fact we need to recall here is
that $\Coder^*(\L^c(\desusp L))$ is a natural pre-Lie algebra, with
the product $\star$ defined as
follows~\cite[Section~II.3.9]{markl-shnider-stasheff:book}. Let
$\Theta,\Omega \in \Coder^*(\L^c(\desusp L))$ and denote by
$\overline{\Omega} : \L^c(\desusp L) \to \desusp L$ the corestriction
of $\Omega$. The pre-Lie product $\Theta \star \Omega$ is then defined
as the coextension of the composition
$$
(-1)^{|\Theta||\Omega|}\cdot
\overline{\Omega}\circ \Theta : \L^c(\desusp L) \to \desusp L,
$$
see~\cite[Section~II.3.9]{markl-shnider-stasheff:book} for details.

By the freeness of the pre-Lie algebra $\pL^*(\bfk^n,\oh)$,
each choice $\Rada f1n \in \Coder(\L^c(\desusp L))$ determines a
unique pre-Lie algebra homomorphism
$$
\Psi_{\Rada f1n} : \pL^*(\bfk^n,\oh) \to \Coder^*(\L^c(\desusp L))
$$  
such that $\Psi_{\Rada f1n}(e_i) := f_i$ for each $1 \leq i \leq n$,
and $\Psi_{\Rada f1n}(\oh) := \lambda$.  Because $\Psi_{\Rada f1n}(\oh
\star \oh)= \lambda^2 = 0$, the map $\Psi_{\Rada f1n}$ induces a map
of the quotient $\rpL^*(\bfk^n) =\pL^*(\bfk^n,\oh)/(\oh \star \oh)$ 
$$
{\rm r}\Psi_{\Rada f1n} : \rpL^*(\bfk^n) \to \Coder^*(\L^c(\desusp L))
$$
Define finally $U(\Rada f1n) \in \CE {d + |f_1| + \cdots +|f_n|}$ by
\begin{equation}
\label{20}
U(\Rada f1n) := {\rm r}\Psi_{\Rada f1n}(U).
\end{equation}
One can easily verify the following formula that relates the
Chevalley-Eilenberg differential $\dce$ with the differential $d$ in
$\rpL^*(n)$:
\begin{eqnarray*}
\lefteqn{
d(U)(f_1,\ldots,f_n) = \dce( U(f_1,\ldots,f_n))} \rule{1em}{0em}
\\
&& \rule{1em}{0em} -
(-1)^{|U|} \sum_{1 \leq i \leq n} 
(-1)^{|f_1| + \cdots + |f_{i-1}|}\cdot
U(f_1,\ldots,\dce( f_i),\ldots, f_n),
\end{eqnarray*}
for each $U \in \rpL^*(n)$ and $\Rada f1n \in \CE*$.

\begin{proposition}
\label{to_jsem_zvedav_jestli_za_Jitkou_pojedu}
The collection 
$\rpL^* := \{\rpL^*(n)\}_{n \geq 1}$ forms an operad in the category of
dg-vector spaces. Formula~(\ref{20}) determines
an action that makes $\CE *$ a differential graded 
$\rpL^*$-algebra. Consequently, the
cohomology operad $H^*(\rpL)$ naturally acts on the Chevalley-Eilenberg
cohomology $H^*_{CE}(L;L)$ of an arbitrary Lie or $L_\infty$ algebra. 
\end{proposition}

\begin{Proof}
The symmetric group $\Sigma_n$ acts on $\rpL^*(n)$ by
permuting the basis $\Rada e1n$ of $\bfk^n$.  The operadic
composition, induced by the vertex
insertion of decorated trees representing elements of $\pL^*(\bfk^n)$, is
constructed by exactly the same method as the one used in the proof
of~\cite[Proposition~II.1.27]{markl-shnider-stasheff:book}. 
The verification that $U$ defines an operadic action is easy.  
\end{Proof}

Let $\Blie^*$ denote, as in Section~1, the dg-operad of
natural operations on the Chevalley-Eilenberg complex of a Lie algebra
with coefficients in itself and $\Bshl^*$ an analog of this operad for
$L_\infty$-algebras. Because each Lie algebra is also an
$L_\infty$-algebra, there exists an obvious
`forgetful' homomorphism $c: \Bshl^* \to \Blie^*$.  By
Proposition~\ref{to_jsem_zvedav_jestli_za_Jitkou_pojedu},
formula~(\ref{20}) defines maps $t : \rpL^* \to \Blie^*$ and
$\widetilde t : \rpL^* \to \Bshl^*$.  The diagram
$$
{
\unitlength=.07em
\begin{picture}(130.00,78.00)(0.00,0.00)
\put(70.00,0.00){\makebox(0.00,0.00){$t$}}
\put(58.00,52.00){\makebox(0.00,0.00){$\widetilde t$}}
\put(130.00,40.00){\makebox(0.00,0.00)[l]{$c$}}
\put(0.00,10.00){\makebox(0.00,0.00){$\rpL^*$}}
\put(120.00,10.00){\makebox(0.00,0.00){$\Blie^*$}}
\put(120.00,70.00){\makebox(0.00,0.00){$\Bshl^*$}}
\put(20.00,10.00){\vector(1,0){80.00}}
\put(120.00,60.00){\vector(0,-1){40.00}}
\put(20.00,20.00){\vector(2,1){80.00}}
\end{picture}}
$$
is clearly commutative and $\widetilde t : \rpL^*
\to \Bshl^*$ is in fact an {\em inclusion\/} of dg-operads, compare the
remarks in Subsection~\ref{mot}.

\section{Proof of Theorem~\protect\ref{.}} 
\label{zitra_se_mozna_poleti}

For the purposes of the proof of Theorem~\ref{.}, it will be
convenient to reduce the complex $\rpL^*(V) = (\rpL^*(V),d)$
constructed in Proposition~\ref{a} as follows. Since the construction
of $\rpL^*(V)$
is functorial in $V$, one may consider the map $\rpL^*(V) \to
\rpL^*(0)$ induced by the map $V \to 0$ from $V$ to the trivial vector
space $0$. The kernel $\overline{\rpL}^*(V) := \Ker(\rpL^*(V) \to
\rpL^*(0))$ is clearly a subcomplex of $\rpL^*(V)$. Since 
$$
\rpL^n(0) = \cases {\Span_\bfk(\oh)}{for $n = 1$ and}0{otherwise,}
$$
the complexes $\overline{\rpL}^*(V)$ and $\rpL^*(V)$ differ only at
the second term, and, under the
isomorphism~(\ref{zitra_odjizdim_na_zavody-uz_tam_sem}),   
$$
\overline{\rpL}^1(V)  \cong  \bigoplus_{n \geq 1} \Tr_n^1(V).
$$
It is also obvious that
$$
\Ker \left(\rule{0em}{1em} d : \pL(V) \to \pL^1(V)  \right)
=  \Ker \left(\rule{0em}{1em} d : \pL(V) \to \overline{\rpL}^1(V)  \right).
$$

The central object of this section is the commutative
diagram:
\begin{equation}
\label{d}
\unitlength .8em
\thicklines
\begin{picture}(21,7)(4,17)
\put(8,17){\makebox(0,0)[cc]{$\pL(V)$}}
\put(8,11){\makebox(0,0)[cc]{$\bT(V)$}}
\put(20,17){\makebox(0,0)[cc]{$\overline{\rpL}^1(V)$}}
\put(20,11){\makebox(0,0)[cc]{$\bT(V) \otimes \bT(V)$}}
\put(8,23){\makebox(0,0)[cc]{$\L(V)$}}
\put(8,16){\vector(0,-1){4}}
\put(20,16){\vector(0,-1){4}}
\put(0,-1){
\put(1,0){
\put(-.5,0){\bezier{50}(8,22)(8,22.25)(7.75,22.25)}
\bezier{50}(7.25,22.25)(7,22.25)(7,22)
}
\put(8,22){\vector(0,-1){3}}
}
\put(9.8,11){\vector(1,0){6.1}}
\put(10,17){\vector(1,0){7.5}}
\put(7,20){\makebox(0,0)[cc]{$i$}}
\put(7,14){\makebox(0,0)[cr]{$p$}}
\put(14,18){\makebox(0,0)[cc]{$d$}}
\put(14,12){\makebox(0,0)[cc]{$\bDelta$}}
\put(21,14){\makebox(0,0)[lc]{$p^1$}}
\end{picture}
\end{equation}
\vskip 5em
\noindent 
in which $i :\L(V) \hookrightarrow \pL(V)$ is the inclusion and $p :
\pL(V) \to \bT(V)_{\it pL} = \bT(V)$ the canonical map of pre-Lie
algebras induced by the inclusion $V \hookrightarrow \bT(V)$. The
definition of $p^1 : \overline{\rpL}^1(V) \to \bT(V) \ot
\bT(V)$ will use the following simple facts.

\vskip .5em
{\em Fact~1\/}. The graded pre-Lie algebra structure of $\rpL^*(V)$
induces on $\rpL^1(V)= \pL^1(V)$ a~structure of
a $\pL(V)$-bimodule. 

{\em Fact~2\/}. With the structure above, $\rpL^1(V)$ is the free
$\pL(V)$-bimodule generated by the dummy~$\oh$.

{\em Fact~3\/}. The $*$-action~(\ref{2}) makes $\T(V) \ot \T(V)$ a
bimodule over the pre-Lie algebra $\T(V)_{pL}$. Therefore $\T(V) \ot
\T(V)$ is a $\pL(V)$-bimodule, via the canonical map $p : \pL(V) \to
\T(V)_{pL}$.

\vskip .5em
The above facts imply that one can define a map $\widehat p^1 : \rpL(V)^1
\to \T(V) \ot \T(V)$ by requiring that it is a $\pL(V)$-bimodule
homomorphism satisfying
$$
\widehat p^1(\oh) := 1 \ot 1 \in \T(V) \ot \T(V). 
$$
It is clear that this $\widehat p^1$ restricts to the requisite map
$p^1 : \brpL^1(V) \to \bT(V) \ot \bT(V)$.

To prove that the bottom square of~(\ref{d}) commutes, we notice that
both compositions $\bDelta p$ and $p^1 d$ behave in the
same way with respect to the pre-Lie multiplication $\star$ on
$\pL(V)$. 
Indeed, for $a,b \in \pL(V))$, by~(\ref{4})
$$
\bDelta p(a \star b) = \bDelta(p(a) \ccdot p(b)) = \bDelta p(a) * p(b)
+ p(a)* \bDelta p(b) + R(p(a),p(b)).
$$
Similarly, by~(\ref{0}) and the definition of $p^1$,
\begin{eqnarray*}
p^1  d(a \star b) &=& p^1(d(a) \star b) + p^1(a \star d(b)) +
p^1(Q(a,b))
\\
&=& p^1 d(a)*  p(b) + p(a) * p^1 d(b) + p^1(Q(a,b)).
\end{eqnarray*}
It remains to verify that $p^1(Q(a,b)) = R(p(a),p(b))$. By the definitions
of $p^1$, $Q$, $R$ and the $*$-action~(\ref{2}),
\begin{eqnarray*}
p^1(Q(a,b)) &=& p^1((\oh \star a) \star b) - p^1(\oh \star (a \star
b)) = p^1(\oh \star a) * p(b) - p^1(\oh) * p(a \star b)
\\
&=&
(p^1(\oh) * p(a)) * p(b) - p^1(\oh) * (p(a) \ccdot p(b))
\\
&=&
((1 \ot 1) * p(a)) * p(b) - (1\ot 1) * (p(a) \ccdot p(b))
\\
&=&
p(a) \ot p(b) + p(b) \ot p(a) = R(p(a),p(b)).
\end{eqnarray*}
Observe finally that $\bDelta p(v) = p^1 dv = 0$ for $v \in V$. The
commutativity $\bDelta p = p^1 d$ of the bottom square
of~(\ref{d}) then follows from the following lemma.

\begin{lemma}
Let $S : \pL(V) \ot \pL(V) \to \bT(V) \ot \bT(V)$ be a symmetric
linear map such that the expression
\begin{equation}
\label{jsem_zvedav_jestli_budou_Pastviny_fungovat}
S(a,b) * p(c) + S(a \star b,c) - S(a,b\star c),\ a,b,c \in \pL(V),
\end{equation}
is symmetric in $b$ and $c$. Then there exists precisely one linear
map $F : \pL(V) \to \bT(V) \ot \bT(V)$ such that
\begin{itemize}
\item[(i)]
$F(a \star b) = F(a)* p(b) + p(a) * F(b) + S(a,b)$ for each $a,b \in
\pL(V)$, and
\item[(ii)]
$F(v) = 0$ for each $v \in V$.
\end{itemize}
\end{lemma}

\begin{Proof}
The map $F$ is constructed by the induction on the
monomial length of elements of $\pL(V)$, its uniqueness is
obvious. The symmetry of the form
in~(\ref{jsem_zvedav_jestli_budou_Pastviny_fungovat}) in $b$ and $c$
is necessary for the compatibility of the rule~(i) with the
axiom~(\ref{aby_to_nezkoncilo_prusvihem}). 
\end{Proof}

We claim that Proposition~\ref{.} follows from the following

\begin{lemma}
\label{aby_mne_nerozbolel_zub}
In diagram~(\ref{d}),
\begin{itemize}
\item[(i)]
$di = 0$,
\item[(ii)]
$\Ker(d) \cap \Ker(p) = 0$ and
\item[(iii)]
$p(\Ker(d)) \subset \L(V)$.
\end{itemize}
\end{lemma}

Indeed, (i)~implies that $\L(V) \subset \Ker(d)$ while (ii) and (iii)
together imply that $p$ maps $\Ker(d)$ monomorphically to
$\L(V)$. Since all these spaces are graded of finite type and
their maps preserve the gradings, one concludes that $\L(V)
= \Ker(d)$.

\noindent 
{\bf Proof of Lemma~\ref{aby_mne_nerozbolel_zub}.} \hglue 1.8em 
The symmetry of $Q$ in~(\ref{0}) implies that $d$ is a derivation of
the Lie algebra $\pL(V)_L$ associated to $\pL(V)$. This fact, together
with $d(V) = 0$, readily implies that $d$ annihilates Lie elements in
$\pL(V)$, which is~(i).

Our proof of~(ii) relies on the tree language introduced in
Section~4. We will use
the following terminology.
A decorated tree $T \in \Tr(V)$ is {\em linear\/} if all its vertices are of
arity $\leq 1$. Such a tree $T$ is of the form
\begin{equation}
\label{18}
{
\unitlength=1.5pt
\begin{picture}(0.00,27.00)(0.00,30.00)
\thicklines
\put(5.00,0.00){\makebox(0.00,0.00)[l]{$v_{i}$}}
\put(5.00,10.00){\makebox(0.00,0.00)[l]{$v_{i-1}$}}
\put(5.00,30.00){\makebox(0.00,0.00)[l]{$v_3$}}
\put(5.00,40.00){\makebox(0.00,0.00)[l]{$v_2$}}
\put(5.00,50.00){\makebox(0.00,0.00)[l]{$v_1$}}
\put(0.00,20.00){\makebox(0.00,0.00){$\vdots$}}
\put(0.00,0.00){\makebox(0.00,0.00){$\bullet$}}
\put(0.00,10.00){\makebox(0.00,0.00){$\bullet$}}
\put(0.00,30.00){\makebox(0.00,0.00){$\bullet$}}
\put(0.00,40.00){\makebox(0.00,0.00){$\bullet$}}
\put(0.00,50.00){\makebox(0.00,0.00){$\root$}}
\put(0.00,10.00){\line(0,-1){10.00}}
\put(0.00,50.00){\line(0,-1){20.00}}
\end{picture}}
\end{equation}
\vskip 40pt
\noindent  
with some $\Rada v1i \in V$, $i \geq 1$.  Each non-linear tree $T \in \Tr(V)$
necessarily looks as
\begin{equation}
\label{19}
{
\unitlength=1.5pt
\begin{picture}(40.00,40.00)(0.00,0.00)
\thicklines
\put(0,-30){
\thicklines
\put(25.00,20.00){\makebox(0.00,0.00)[l]{$v_{i}$}}
\put(25.00,30.00){\makebox(0.00,0.00)[l]{$v_{i-1}$}}
\put(25.00,50.00){\makebox(0.00,0.00)[l]{$v_3$}}
\put(25.00,60.00){\makebox(0.00,0.00)[l]{$v_2$}}
\put(25.00,70.00){\makebox(0.00,0.00)[l]{$v_1$}}
\put(20.00,10.00){\makebox(0.00,0.00){$S$}}
\put(20.00,40.00){\makebox(0.00,0.00){$\vdots$}}
\put(20.00,20.00){\makebox(0.00,0.00){$\bullet$}}
\put(20.00,30.00){\makebox(0.00,0.00){$\bullet$}}
\put(20.00,50.00){\makebox(0.00,0.00){$\bullet$}}
\put(20.00,60.00){\makebox(0.00,0.00){$\bullet$}}
\put(20.00,70.00){\makebox(0.00,0.00){$\bullet$}}
\put(40.00,0.00){\line(-1,1){20.00}}
\put(0.00,0.00){\line(1,0){40.00}}
\put(20.00,20.00){\line(-1,-1){20.00}}
\put(20.00,30.00){\line(0,-1){10.00}}
\put(20.00,70.00){\line(0,-1){20.00}}
}
\end{picture}}
\end{equation}
\vskip 40pt
\noindent 
where $S$ is a tree whose root vertex $v_i$ has arity $\geq 2$. We say
that such a decorated tree has {\em tail of length $i$\/}. These
notions translate to decorated trees from $\Tr^1(V)$ in the obvious
manner.

We leave to the reader to verify that, under
identification~(\ref{zitra_odjizdim_na_zavody}), the map $p : \pL(V)
\to \bT(V)$ is described as
$$
p(T) =
\cases{v_1 \ot v_2 \ot \cdots \ot v_{i-1} \ot v_i}
      {if $T$ is linear as in~(\ref{18}), and}0
      {if $T$ is non-linear.}
$$
Therefore $\Ker(p)$ consists of linear combinations of non-linear
trees. Before going further, we need to inspect how the map $\delta :
\Tr(V) \to \Tr^1(V)$ of~(\ref{Nunyk}), which is the differential $d:
\pL(V) \to \pL^1(V)$ written in terms of trees, acts on non-linear
trees.  If $T$ is the decorated tree~(\ref{19}), then it immediately
follows from the definition~(\ref{Nunyk}) of $\delta$ that
\begin{equation}
\label{snad-to-konecne-dopisu}
\delta (T) = -T' + \mbox {trees with tails of length $\leq  i$},
\end{equation}
where $T'$ is the following decorated tree with tail of length $i+1$
$$
{
\unitlength=1.5pt
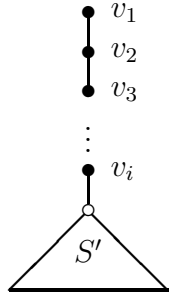
\begin{picture}(40.00,75.00)(0.00,0.00)
\thicklines
\put(25.00,30.00){\makebox(0.00,0.00)[l]{$v_i$}}
\put(25.00,50.00){\makebox(0.00,0.00)[l]{$v_3$}}
\put(25.00,60.00){\makebox(0.00,0.00)[l]{$v_2$}}
\put(25.00,70.00){\makebox(0.00,0.00)[l]{$v_1$}}
\put(20.00,10.00){\makebox(0.00,0.00){$S'$}}
\put(20.00,40.00){\makebox(0.00,0.00){$\vdots$}}
\put(20.00,20.00){\makebox(0.00,0.00){$\oh$}}
\put(20.00,30.00){\makebox(0.00,0.00){$\bullet$}}
\put(20.00,50.00){\makebox(0.00,0.00){$\bullet$}}
\put(20.00,60.00){\makebox(0.00,0.00){$\bullet$}}
\put(20.00,70.00){\makebox(0.00,0.00){$\bullet$}}
\put(40.00,0.00){\line(-1,1){19.00}}
\put(0.00,0.00){\line(1,0){40.00}}
\put(19.00,19.00){\line(-1,-1){19.00}}
\put(20.00,30.80){\line(0,-1){9.00}}
\put(20.00,70.00){\line(0,-1){20.00}}
\end{picture}}
$$
in which $S'$ is the tree obtained from $S$ by replacing the root
vertex decorated by $v_i$ by the special one. 
The map $\delta' : \Tr(V) \to \Tr^1(V)$ given by $\delta'(T) :=
-T'$ is a {\em monomorphism\/}.

Let $x$ be a linear combination of non-linear trees and assume $\delta
(x) =0$. We must prove that then $x=0$. Assume $x \not= 0$ and
decompose $x = x_s + x_{s-1} + \cdots + x_1$, where $x_i$ is, for $1
\leq i \leq s$, a linear combination of decorated trees with tails of
length $i$, and $x_s \not= 0$.  By~(\ref{snad-to-konecne-dopisu}), the
only trees with tails of length $s+1$ in $\delta(x)$ are those spanning
$\delta'(x_s)$, therefore $\delta(x) = 0$ implies $\delta'(x_s)=0$ which
in turn implies that $x_s = 0$, because $\delta'$ is monic. This is a
contradiction, therefore $x=0$ which proves~(ii).

To verify~(iii), notice that, by the commutativity of the bottom
square of~(\ref{19}), $p (\Ker(\delta)) \subset \Ker(\bDelta))$ while
$\Ker(\bDelta) = \L(V)$ by Theorem~\ref{psano_v_aute}. This finishes
the proof of the lemma.
\qed

\section{Some open questions and ramifications}

\begin{odstavec}
\label{trip}
{\rm
{\it Triplettes of operads (after J.-L.~Loday)\/}.
The following notion was introduced in~\cite{loday:slides}. 

\begin{Definition}
\label{jeste_ji_musim_napsat_SMS}
The data $(\calC,\spin,\Aalg \stackrel{F}\to \Palg)$, where
\begin{itemize}
\item[(i)]
$\calC$ and  $\calA$ are operads,
\item[(ii)]
$\spin$ are `spin' relations intertwining $\calC$-co-operations and
$\calA$-operations, so that $(\calC,\spin,\calA)$ determines a class of
bialgebras, 
\item[(iii)]
the operad $\calP$ governs the algebra structure of the primitive
part $\Prim(\calH)$ of  $(\calC,\spin,\calA)$-bialgebras, and
\item[(iv)]
$F$ is a forgetful functor functor from the category of
$\calA$-algebras to the category of $\calP$-algebras such that the
inclusion $\Prim(\calH) \subset F(\calH)$ is a morphism of $\calP$-algebras,
\end{itemize}
is called a {\em triplette\/} of operads.
\end{Definition}

An example is $(\Com,\spin, \Ass,\Lie)$, with $\spin$
the usual bialgebra relation recalled
in~(\ref{to_jsem_zvedav_jesli_budou_pastviny_fungovat}). 
Let $U$ be a left adjoint to $F$.
A triplette in Definition~\ref{jeste_ji_musim_napsat_SMS} is 
{\em good\/}~\cite{loday:slides}, 
if the following three conditions are equivalent:

\begin{itemize}
\item[(i)]
a $(\calC,\spin,\calA)$-bialgebra $\calH$ is connected,
\item[(ii)]
$\calH \cong U(\Prim(\calH))$, and
\item[(iii)]
$\calH$ is cofree among connected
$\calC$-coalgebras.
\end{itemize}
Let $\calA(V)$ (resp.~$\calP(V)$) denote the free $\calA$-
(resp.~$\calP$-)algebra on $V$.
As observed in~\cite{loday:slides}, for good triplettes
\begin{equation}
\label{`}
\Prim (\calA (V)) \cong \calP(V).
\end{equation}
 
The classical Theorem~\ref{psano_v_aute} in Section~2 is a
consequence of the goodness of the triplette $(\Com,\spin, \Ass,\Lie)$
mentioned above, because~(\ref{`}) in this case says that $\Prim(T(V))
\cong \L(V)$. Other, in some cases very surprising, good triplettes
can be found in~\cite{loday:slides}.  The following problem was
suggested by J.-L.~Loday:

\begin{problem}
\label{kk}
Are there an
operad $\calC$ and spin relations $\spin$ with the property that
$(\calC,\spin,\pLie,\Lie)$ is a good triplette?
\end{problem}

As we remarked in Subsection~\ref{gen}, the affirmative answer to the
Deligne conjecture given in~\cite{kontsevich-soibelman} implies that
there exist a characterization of Lie elements in brace
algebras~\cite{gerstenhaber-voronov:FAP95} similar to our
Theorem~\ref{.}. This suggests formulating the following version of
Problem~\ref{kk} in which $\Brace$ is the operad for brace
algebras.

\begin{problem}
Are there an
operad $\calC$ and spin relations $\spin$ with the property that
$(\calC,\spin,\Brace,\Lie)$ is a~good triplette?
\end{problem}

}\end{odstavec}

\begin{odstavec}
{\rm 
{\it Lie elements and cobar constructions.\/} 
In Section~2 we calculated the cohomology
of the cobar construction~(\ref{1a}) of the
shuffle coalgebra and observed that $H^0(\T(V),\Delta)$ is
isomorphic to the free Lie algebra $\L(V)$. In our characterization of
Lie elements in pre-Lie algebras, the role of~(\ref{1a}) is played by
complex~(\ref{za_chvili_tam_musim_volat_a_hrozne_se_mi_nechce!}).
This leads to the following problem, which may or may not be related
to Problem~\ref{kk},

\begin{problem}
Calculate the cohomology
of~(\ref{za_chvili_tam_musim_volat_a_hrozne_se_mi_nechce!}). Is 
this complex the cobar construction of some coalgebra?
\end{problem}

As D.~Tamarkin recently informed us, methods proposed in an enlarged
unfinished, 
unpublished version of~\cite{tamarkin:def_of_chiral_algebras} may
imply that the
complex~(\ref{za_chvili_tam_musim_volat_a_hrozne_se_mi_nechce!}) is
acyclic in positive dimensions, as envisaged also by some
conjectures formulated in~\cite{markl:de}.  }\end{odstavec}

\references

\nextref{chapoton-livernet:pre-lie}
        {Chapoton, F. and Livernet, M.}
        {\em Pre-{Lie} algebras and the rooted trees operad}
        {Internat. Math. Res. Notices {\bf 8}(2001), 395--408}

\nextref{chevalley-eilenberg}
        {Chevalley, C. and Eilenberg, S.}
        {\em Cohomology theory of {Lie} groups and {Lie} algebras}
        {Trans. Amer. Math. Soc., {\bf 63}(1948), 85--124}

\nextref{dzhu-lof:HHA02}
        {Dzhumadil'daev, A. and L\"ofwall, C.}
        {\em Trees, free right-symmetric algebras, free {Novikov} algebras and
         identities}
        {Homotopy, Homology and Applications, {\bf 4(2)}(2002),
         165--190}

\nextref{gerstenhaber:AM63}
        {Gerstenhaber, M.}
        {\em The cohomology structure of an associative ring}
        {Ann. of Math. {\bf 78(2)}(1963), 267--288}

\nextref{gerstenhaber-voronov:FAP95}
        {Gerstenhaber, M. and Voronov, A.A.}
        {\em Higher operations on the {Hochschild} complex}
        {Functional Anal. Appl. {\bf 29(1)}(1995), 1--6 (in Russian)}

\nextref{guin-oudom}
        {Guin, D. and Oudom, J.-M.}
        {\em Sur l'alg\'ebre enveloppante d'une alg\'ebre pr\'e-{Lie}}
        {C. R. Acad. Sci. Paris S{\'e}r. I Math., {\bf 340(5)}{2005}, 331--336}

\nextref{kontsevich-soibelman}
        {Kontsevich, M. and Soibelman, Y.}
        {\em Deformations of algebras over operads the {Deligne} conjecture}
        {In {Dito, G. et al.}, editor, {Conf\'erence {Mosh\'e Flato} 1999:
         Quantization, deformation, and symmetries}, number~21 in 
         Math. Phys. Stud.,
         pages 255--307. Kluwer Academic Publishers, 2000}

\nextref{lada-markl:CommAlg95}
        {Lada, T. and Markl, M.}
        {\em Strongly homotopy {Lie} algebras}
        {Comm. Algebra {\bf 23(6)}(1995), 2147--2161}

\nextref{lm:sb}
        {Lada, T. and Markl, M.}
        {\em Symmetric brace algebras}
        {Applied Categorical Structures {\bf 13(4)}, 351--370}

\nextref{loday:slides}
        {Loday, J.-L.}
        {\em Generalized bialgebras and triples of operads}
        {A slide show. Available from the {\tt www} home page of
         J.-L.~Loday}

\nextref{maclane:homology}
        {Mac~Lane, S.}
        {\em Homology}
        {Springer-Verlag, 1963}

\nextref{markl:de}
        {Markl, M.}
        {\em Cohomology operations and the {Deligne} conjecture}
        {preprint {\tt math.AT/0506170}, June 2005}

\nextref{markl-remm}
        {Markl, M. and Remm, E.}
        {Algebras with one operation including {Poisson} and other
         {Lie}-admissible algebras}
        {J. Algebra, {\bf 299}(2006), 171--189}

\nextref{markl-shnider-stasheff:book}
        {Markl, M., Shnider, S. and Stasheff, J.D.}
        {\em Operads in Algebra, Topology and Physics}
        {volume~96 of Mathematical Surveys and Monographs,
         American Mathematical Society, Providence, Rhode Island,
         2002}

\nextref{ree:AnM69}
        {Ree, R.}
        {\em Lie elements and an algebra associated with shuffles}
        {Ann. of Math. {\bf 68}, 210--220}

\nextref{serre:65}
        {Serre, J.-P.}
        {\em Lie Algebras and {Lie} Groups}
        {Benjamin, 1965}

\nextref{tamarkin:def_of_chiral_algebras}
        {Tamarkin, D.}
        {\em Deformations of chiral algebras}
        {In {Ta Tsien et al.}, editor, {Proceedings of ICM 2002, Beijing, 
         China, August 20-28},
         vol. II: Invited lectures, pages 105--116. 
         Beijing: Higher Education Press;
         Singapore: World Scientific/distributor, 2002}

\lastpage
\end{document}